\documentclass[10pt,onecolumn,twoside]{IEEEtran}

\pdfoutput=1
\usepackage{amsmath, amssymb, amsfonts}  
\usepackage{cases}

\usepackage{graphicx}    

\usepackage{bm}
\usepackage{bbm}
\usepackage{dsfont}
\usepackage{xcolor}

\usepackage{ifpdf}

\usepackage{enumitem}  
\setlist[description]{leftmargin=\parindent}

\usepackage[mathscr]{eucal}
\usepackage[caption=false]{subfig}
\usepackage{theorem}

\newtheorem{theorem}{Theorem}

\newtheorem{proposition}{Proposition}
\newtheorem{lemma}{Lemma}
\newtheorem{corollary}{Corollary}

{\theorembodyfont{\upshape}
\newtheorem{remark}{Remark}
\newtheorem{example}{Example}

}

\def\qed{\hfill \vrule height 5pt width 5pt depth 0pt \medskip}
\newcommand{\proof}{\noindent {\bf Proof. }}

\newcommand{\beq}{\begin{equation}}
\newcommand{\eeq}{\end{equation}}
\newcommand{\beqa}{\begin{eqnarray}}
\newcommand{\eeqa}{\end{eqnarray}}
\newcommand{\beqan}{\begin{eqnarray*}}
\newcommand{\eeqan}{\end{eqnarray*}}

\newcommand{\pde}[2]{ \frac{\partial #1}{\partial #2} }
\newcommand{\ppde}[2]{\frac{\partial^2 #1}{\partial #2^2}}

\newcommand{\abs}[1]{ \lvert #1 \rvert }
\newcommand{\norm}[1]{ \| #1 \|}

\newenvironment{system}{\left\lbrace\begin{array}{@{}l@{}}}{\end{array}\right.}

\newcommand{\R}{\mathbb{R}}
\newcommand{\C}{\mathbb{C}}
\newcommand{\K}{\mathcal{K}}
\newcommand{\1}{\mathbbm{1}}

\newcommand{\diag}{\text{diag}}

\begin{document}

\title{Multiequilibria analysis for a class of collective decision-making networked systems%
\thanks{Work supported in part by a grant from the Swedish Research Council (grant n. 2015-04390).}}

\author{Angela~Fontan~and~Claudio~Altafini
	\thanks{A. Fontan and C. Altafini are with the Division of Automatic Control,
		Department of Electrical Engineering,
		Link\"{o}ping University,
		SE-58183 Link\"{o}ping,
		Sweden,
		E-mail: $\{$angela.fontan, claudio.altafini$\}$@liu.se }}

\maketitle

\begin{abstract}
The models of collective decision-making considered in this paper are nonlinear interconnected cooperative systems with saturating interactions.
These systems encode the possible outcomes of a decision process into different steady states of the dynamics.
In particular, they are characterized by two main attractors in the positive and negative orthant, representing two choices of agreement among the agents, associated to the Perron-Frobenius eigenvector of the system.
In this paper we give conditions for the appearance of other equilibria of mixed sign. 
The conditions are inspired by Perron-Frobenius theory and are related to the algebraic connectivity of the network.
We also show how all these equilibria must be contained in a solid disk of radius given by the norm of the equilibrium point which is located in the positive orthant.
\end{abstract}


\section{Introduction}
\IEEEPARstart{N}{onlinear} interconnected systems are used in broadly different contexts to describe the collective dynamical behavior of an ensemble of ``agents'' interacting with each other in a non-centralized manner.
They are used for instance to represent collective decision-making by animal groups \cite{bib:FranciSrivastavaLeonard2015,Gray2017Agent,Leonard14Multi}, formation of opinion in social communities \cite{bib:AltafiniLini2015,bib:Altafini2012}, dynamics of a gene regulatory network \cite{bib:LiChenAihara2006}, or neural networks \cite{bib:Haykin1999}.
Such models often share similar features, like the fact of using first order dynamics at a node and sigmoidal saturation-like functions to describe the interactions among the nodes \cite{bib:Hopfield1984,bib:Haykin1999,Kaszkurewicz2000Matrix,siljak1978large}.
The latter functional form is instrumental to avoid diverging dynamics. 
The price to pay for having an effectively ``bounded'' dynamics is however the appearance of complex dynamical phenomena such as periodic orbits or multiple equilibrium points, which complicate considerably the behavior of the system and its understanding.
While (stable) periodic orbits can be ruled out by choosing functional forms which are, beside saturated, also monotone \cite{bib:HirschSmith2005,bib:Smith1988}, it is in general more difficult to deal with multiple equilibria.
These are often a necessary feature of a model, like for instance, when describing bistability in a biological system. 
In the literature on agent-based animal groups for instance, the collective decision is typically a selection between two different attractors \cite{bib:FranciSrivastavaLeonard2015,Gray2017Agent,Leonard14Multi}. 
A scenario that occurs often is that of a finite number of attractors, each with its own basin of attraction, encoding the possible outcomes of the collective decision-making process. 
This is for instance the model of choice of several classes of neural networks, like the so-called Hopfield \cite{bib:Hopfield1984} or Cohen-Grossberg \cite{bib:CohenGrossberg1988} neural networks. 
When a neural network is interpreted as an associative memory storage device or is used for pattern recognition, the presence of a high number of stable equilibria increases the storage capacity of the network.
The multistability of neural networks has in fact been extensively investigated in recent years \cite{bib:ChengLinShih2006,bib:ChengLinShih2007,bib:ZengZheng2012,bib:LuWangChen2011,bib:FortiTesi2015}, see also \cite{bib:ZhangWangLiu2014} for an overview.
%

The model adopted in this paper is described in \cite{bib:FranciSrivastavaLeonard2015,UgoAbara2016Spectral,UgoAbara2015Existence}.
All interactions are ``activatory'', i.e., the adjacency matrix is nonnegative.
In addition, it is symmetric or diagonally symmetrizable and irreducible. 
The model has a Laplacian-like structure at the origin and monotone saturating nonlinearities to represent the interaction terms. 
The amplitude of the interaction part is modulated by a scalar parameter, interpretable as the strength of the social commitment of the agents, and playing the role of bifurcation parameter.
When the parameter is small, the origin is globally stable, as can be easily deduced by (global) diagonal dominance. 
The interesting dynamics happens when the parameter passes a bifurcation threshold: the origin becomes unstable, and two locally stable equilibria, one positive, the other negative, are created.
This is the behavior described in \cite{bib:FranciSrivastavaLeonard2015}.
The bifurcation analysis of \cite{bib:FranciSrivastavaLeonard2015}, however, captures what happens only in a neighborhood of the bifurcation point. 
Moving away from the bifurcation, all is known is that the positive and negative orthants remain invariant sets for the dynamics, and each keeps having a single asymptotically stable equilibrium point, see \cite{UgoAbara2016Spectral,UgoAbara2015Existence}. 
What happens in the remaining orthants is unknown and its investigation is the scope of this paper. 

It is useful to look at the neural network literature (in our knowledge the only field that has studied the multiequilibria problem is a systematic way for such interconnected systems). 
It is known since \cite{bib:Hopfield1984} that for neural networks with connections that are symmetric, monotone increasing and sigmoidal, the equations of motion always lead to convergence to stable steady states.
The number of such equilibria is shown to grow exponentially with the number $n$ of ``neurons'' for various specific models \cite{bib:ChengLinShih2006,bib:ChengLinShih2007,bib:ZengZheng2012,bib:LuWangChen2011}.
In order to count equilibria, often in these papers one has to resort to nondecreasing piecewise linear functions to describe the saturations, and to obtain algebraic conditions on the equilibria by looking at the corners and at the constant slopes. 
The existence of the equilibria is also checked through Brouwer fixed-point arguments.
None of these methods apply in our case. 

The consequence is that in order to investigate the presence of multiple equilibria in our system we have to take a completely different approach. 
Since the adjacency matrix of our network is nonnegative, we can use Perron-Frobenius theorem and the geometrical considerations that follow from it.
The main result of the paper is a necessary and sufficient condition for existence of equilibria outside $ \mathbb{R}^n_+/\mathbb{R}^n_-$, formulated in terms of the second largest positive eigenvalue of the adjacency matrix. 
Roughly speaking, this condition says that the interval of values of the bifurcation parameter in which no mixed-sign equilibrium can appear is determined by the algebraic connectivity of the Laplacian of the system \cite{Deabreu2007Old,Donetti2006Optimal,Olfati2003Consensus}. 
For consensus problems, the role of the algebraic connectivity is well-known: the bigger is the gap between $0$ (the least eigenvalue of the Laplacian) and the algebraic connectivity, the more robust the consensus is to model uncertainties, parameter variations, node or link failures, etc. \cite{Donetti2006Optimal}.
In the present context, the spectral gap in the Laplacian (or, more properly, in the adjacency matrix), represents the range of values of social commitment of the agents which leads to a choice between two alternative agreement solutions, with a guaranteed global convergence. Beyond the value represented by the algebraic connectivity, however, the system bifurcates again, and a number of (stable and unstable) mixed-sign equilibria appears quickly, which destroys the global convergence to the agreement manifold in which the two alternative attractors live.  
Although we can investigate these extra equilibria only numerically, it is nevertheless possible to compute exactly the region in which they must be.
The second analytical result of the paper is in fact that for all values of the bifurcation parameter the mixed-sign equilibria have to have norm less than the norm of the equilibrium in $ \mathbb{R}^n_+$. 
In other words, all the equilibria of the system must be contained in a ball centered in the origin of $ \R^n $ and of radius equal to the norm of the positive equilibrium point. 
Our numerical analysis shows that stable equilibria tend to localize towards the boundary of this ball, while those near the origin tend to have Jacobians with several positive eigenvalues and a similar number of positive and negative components.

The rest of this work is organized as follows.  Preliminary material is introduced in Section \ref{section:Preliminaries}. The main results are stated in Sections~\ref{section:MultipleEquilibriaNonlinearCooperativeSystems} and~\ref{section:EquilibriumPointPositiveOrthant}; the first section provides the necessary and sufficient condition for the existence of mixed-sign equilibria, while the second describes the region in which they must be contained. 
Section \ref{section:NumericalAnalysis} finally provides a numerical analysis of the equilibria. 

\section{Preliminaries}
\label{section:Preliminaries}

\subsection{Concave and convex functions}
\label{section:ConcaveConvexFunctions}
Let $\mathcal U\in \R^n$ be a convex set. A function $f:\mathcal U\to \R$ is {\itshape convex} if for all $x_1,x_2\in \mathcal U$ and $\theta$ with $0\le \theta\le 1$, we have
\begin{equation}
	f(\theta x_1+(1-\theta)x_2)\le \theta f(x_1)+(1-\theta)f(x_2).
\label{def:convexFunction}
\end{equation}
It is {\itshape strictly convex} if strictly inequality holds in \eqref{def:convexFunction} whenever $x_1\ne x_2$ and $0<\theta<1$. We say $f$ is {\itshape concave} if $-f$ is convex, and {\itshape strictly concave} if $-f$ is strictly convex.
Suppose $f:\mathcal U\to \R$ is differentiable. Then $f$ is convex if and only if $\mathcal U$ is convex and
\begin{equation}
	f(x_1)\ge f(x_2)+\frac{\partial f}{\partial x}(x_2) (x_1-x_2)
\label{eqn:convexityFirstOrder}
\end{equation}
for all $x_1,x_2\in\mathcal U$. It is strictly convex if and only if $\mathcal U$ is a convex set and strict inequality holds in \eqref{eqn:convexityFirstOrder} for all $x_1\ne x_2$ and $x_1,x_2\in \mathcal U$.


\subsection{Nonnegative matrices and Perron-Frobenius}
\label{subsection:LocationEigenvaluesNonnegativeMatrices}
The set of all $\lambda \in \C$ that are eigenvalues of $A\in \R^{n \times n}$ is called the {\itshape spectrum} of $A$ and is denoted by $\Lambda(A)$.
The {\itshape spectral radius} of $A$ is the nonnegative real number $\rho(A)=\max\{|\lambda|:\,\lambda \in \Lambda(A)\}$.

A matrix $B\in \R^{n\times n}$ is said to be {\itshape similar} to a matrix $A\in \R^{n\times n}$, abbreviated $B\sim A$, if there exists a nonsingular matrix $S\in \R^{n\times n}$ such that $B=S^{-1}AS$.
If $A$ and $B$ are similar, then they have the same eigenvalues, counting multiplicity.

A matrix $A \in \R^{n \times n}$ is said to be {\itshape reducible} if either $n=1$ and $A=0$ or if $n\ge 2$, there is a permutation matrix $P \in \R^{n \times n}$ and there is some integer $r$ with $1\le r\le n-1$ such that $P^TAP=\begin{bmatrix}
	B & C \\ 0 & D
	\end{bmatrix}$
where $B \in \R^{r \times r}$, $C \in \R^{r \times (n-r)}$, $D \in \R^{(n-r) \times (n-r)}$. 
A matrix $A \in \R^{n \times n}$ is said to be {\itshape irreducible} if it is not reducible.

The matrix $A=[a_{ij}]\in \R^{n \times n}$ is said to be {\itshape diagonally dominant} if $|a_{ii}|\ge \sum_{j\ne i} |a_{ij}|$ for all $i$.
It is said to be {\itshape strictly diagonally dominant} if $|a_{ii}|> \sum_{j\ne i} |a_{ij}|$ for all $i$.

\begin{theorem}[\cite{bib:HornJohnson1990}, Theorem 6.1.10]
	\label{theorem:SDDmatrix}
	Let $A \in \R^{n \times n}$ be strictly diagonally dominant. Then $A$ is nonsingular. If $a_{ii}>0$ for all $i=1,\dots,n$, then every eigenvalue of $A$ has positive real part. If $A$ is symmetric and $a_{ii}>0$ for all $i=1,\dots,n$, then $A$ is positive definite.
\end{theorem}

\begin{theorem}[Perron-Frobenius, \cite{bib:BermanPlemmons1994} Theorem 1.4]
	\label{theorem:Perron-Frobenius}
	If $A\in \R^{n \times n}$ is irreducible and nonnegative then $\rho(A)$ is a real, positive, algebraically simple eigenvalue of $A$, of right (left) eigenvector $v>0$ ($w>0$).
	Furthermore, for every eigenvalue $\lambda\in \Lambda(A)$ such that $\lambda\ne \rho(A)$ it is $|\lambda|<\rho(A)$ and the corresponding eigenvector $ v_\lambda $ cannot be nonnegative.
\end{theorem}
\begin{corollary}[Perron-Frobenius]
	\label{corollary:Perron-Frobenius}
	If $A\in \R^{n \times n}$ is irreducible and nonnegative then either
	\begin{equation}
	\label{eq:PF_equal1}
	\rho(A)=\sum_{j=1}^n a_{ij},\quad i=1,\dots,n
	\end{equation}
	or
	\begin{equation}
	\label{eq:PF_disequal1}
	\min_i \left(\sum_{j=1}^n a_{ij}\right)
	\le \rho(A)
	\le \max_{i} \left(\sum_{j=1}^n a_{ij}\right).
	\end{equation}
\end{corollary}


\subsection{Symmetric, symmetrizable and congruent matrices}
\label{subsection:SymmetricSymmetrizableCongruentMatrices} 
Let $A\in \R^{n\times n}$ be symmetric. Then all the eigenvalues of $A$ are real and $S^TAS$ is symmetric for all $S\in \R^{n\times n}$.

A matrix $A\in \R^{n\times n}$ is {\itshape (diagonally) symmetrizable} if $DA$ is symmetric for some diagonal matrix $D$ with positive diagonal entries \cite{bib:Gruber1995}.
The matrix $DA$ is called {\itshape symmetrization} of $A$ and the matrix $D$ is called {\itshape symmetrizer} of $A$ \cite{bib:DiasCastonguayDourado2016}.
The eigenvalues of a symmetrizable matrix are real.
$A\in \R^{n\times n}$ is symmetrizable if and only if it is sign symmetric, i.e. $a_{ij}=a_{ji}=0$ or $a_{ij}\cdot a_{ji}>0$, $\forall \;i\ne j$, and
$a_{i_1,i_2}a_{i_2,i_3}\dots a_{i_k,i_1}=a_{i_2,i_1}a_{i_3,i_2}\dots a_{i_1,i_k}$
for all $i_i,\dots,i_k$.


A matrix $B\in \R^{n\times n}$ is said to be {\itshape congruent} to the matrix $A\in \R^{n\times n}$ if there exists a nonsingular matrix $S$ such that $B=SAS^T$. 
\begin{theorem}[Ostrowski, Theorem 4.5.9 in \cite{bib:HornJohnson1990}]
	\label{theorem:Ostrowski}
	Let $A,S\in \R^{n\times n}$ with $A$ symmetric and $S$ nonsingular. Let the eigenvalues of $A$, $SAS^T$ and $SS^T$ be arranged in nondecreasing order. For each $k=1,\dots,n$ there exists a positive real number $\theta_k$ such that $\lambda_1(SS^T)\le \theta_k \le \lambda_n(SS^T)$ and $\lambda_k(SAS^T)=\theta_k\lambda_k(A)$.
\end{theorem}


\subsection{Cooperative systems}
\label{subsection:CooperativeSystem}
Consider the system
\begin{equation}
\dot x=f(x),\quad x(0)=x_0
\label{eqn: System theory}
\end{equation}
where $f$ is a continuously differentiable function defined on a convex, open set $U\subseteq \R^n$.
We write $x(t,x_0)$ for the forward solution of \eqref{eqn: System theory} with initial condition $x_0\in \R^n$ at $t=0$.

The set $\K$ given by
$\K=\{x \in \R^n: k_i x_i\ge 0,\;k_i\in \{\pm 1\},\,i=1,\dots,n\}$
is an orthant in $\R^n$.
$\K$ is a cone in $\R^n$ and it generates a partial ordering ``$\le_{\K}$'', i.e. $x \le_{\K} y $ iff $y-x \in \K$. The subscript ``$\K$'' will be dropped in case $\K=\R^n_+$, the nonnegative orthant. 

System \eqref{eqn: System theory} is said to be type-$\K$ {\itshape monotone} \cite{bib:Smith1988} if whenever $\bar x$ and $\bar y$ lie in $U$ and if $\bar x \le_\K \bar y$ then $x(t,\bar x)\le_\K y(t,\bar y)$ for all $t \ge 0$ for which both solutions are defined.
In this case we say that the flow of \eqref{eqn: System theory} preserves the ordering $\le_\K$.
\begin{lemma}[\cite{bib:Smith1988}]
	If $f \in C^1(U)$ where $U$ is open and convex in $\R^n$, then $x(t,x_0)$ preserves the partial ordering $\le_\K$ for $t\ge 0$ if and only if $K\frac{\partial f}{\partial x}(x) K$ has nonnegative off-diagonal elements for every $x\in U$, where $K=\diag\{k_1,\dots,k_n\}$, $k_i\in\{\pm 1\}$, is the ``signature'' of the orthant $\K$.
\end{lemma}

If $\K=\R^n_+$ then we have the class of cooperative systems.
System \eqref{eqn: System theory} is said to be {\itshape cooperative} in $U \subset \R^n$ if the differentiable vector field $f: \R^n\to \R^n$ is such that the Jacobian matrix $\frac{\partial f}{\partial x}(x)$ is Metzler for all $x \in U$, that is $\left[\frac{\partial f}{\partial x}(x)\right]_{ij}\ge 0$ for all $i\ne j$.


\section{Multiple equilibria in collective decision-making ystems}
\label{section:MultipleEquilibriaNonlinearCooperativeSystems}
The class of nonlinear systems considered in this work is the following \cite{bib:FranciSrivastavaLeonard2015,UgoAbara2016Spectral}
\begin{equation}
\dot x =-\Delta x+\pi A \psi(x)
\label{eqn:SystemSFA}
\end{equation}
where $x\in \R^n$, $\pi>0$ is a scalar parameter, $A=[a_{ij}]$ is the weighted adjacency matrix of the network, $\Delta=\diag\{\delta_1,\dots,\delta_n\}$, and $\psi(x)=[\psi_1(x_1)\,\dots\, \psi_n(x_n)]^T$.
The matrix $A$ is assumed to be {\itshape nonnegative} with {\itshape null diagonal}, {\itshape irreducible} and {\itshape symmetrizable}.
A Laplacian-like assumption links $\Delta $ and $ A$: $\delta_i=\sum_{j} a_{ij}$.
In the context of agent-based group decisions, $ \delta_i $ represents the inertia of the $i$-th agent to the development of an opinion, $ \psi_i(x_i) $ the capacity of the $i$-th agent of transmitting its opinion to the other agents, mediated by the pairwise susceptibilities $ a_{ji} $.
The parameter $ \pi $ represents a community social effort. See \cite{bib:FranciSrivastavaLeonard2015,Leonard14Multi} for more details.

The vector of functions $\psi(x)$ is such that each $\psi_i(x_i):\R\to \R$ satisfies the following conditions
\begin{gather}
\psi_i(x_i)=-\psi_i(-x_i),\,\forall x_i\in \R\quad\text{(odd)}
\tag{A.1}
\label{assumption:1psiOdd}
\\
\frac{\partial \psi_i}{\partial x_i}(x_i)>0,\,\forall x_i\in \R\;\text{and }\frac{\partial \psi_i}{\partial x_i}(0)=1\quad\text{(monotone)}
\tag{A.2}
\label{assumption:2psiMonotone}
\\
\lim_{x_i\to\pm \infty} \psi_i(x_i)=\pm 1\quad\text{(saturated)}.
\tag{A.3}
\label{assumption:3psiSaturated}
\end{gather}
The assumption \eqref{assumption:3psiSaturated} guarantees the boundedness of the solutions, and together with \eqref{assumption:2psiMonotone} allows to exclude the presence of limit cycles.
Typical choices for $ \psi_i$ are a hyperbolic tangent function, a (modified) Michaelis-Menten function $\psi_i(x_i)=\frac{x_i}{1+|x_i|},\, x_i \in \R$ \cite{bib:Altafini2012}, or a (modified) Boltzmann function $\psi_i(x_i)=\frac{1-e^{-2x_i}}{1+e^{-2x_i}},\, x_i \in \R$ \cite{bib:Haykin1999}. 
The versions here proposed satisfy the conditions \eqref{assumption:1psiOdd}$\div$\eqref{assumption:3psiSaturated}.

Additionally, if each nonlinear function $\psi_i(x_i)$ satisfies the following condition
\beq
\psi_i(x_i) \; 
\begin{cases}
\text{strictly convex} & \forall\, x_i<0\\
\text{strictly concave}& \forall\, x_i>0
\end{cases} \quad \text{(sigmoidal)}
\tag{A.4}
\label{assumption:4Sigmoidal}
\eeq
then, from \eqref{eqn:convexityFirstOrder}, $\abs{\psi_i(x_i)}< \abs{x_i}$ for all $i$ and $x_i\ne 0$.
While all the examples mentioned above (hyperbolic tangent, Michaelis-Menten and Boltzmann functions) satisfy also \eqref{assumption:4Sigmoidal}, nonsigmoidal functions often satisfy \eqref{assumption:1psiOdd}$\div$\eqref{assumption:3psiSaturated} but not \eqref{assumption:4Sigmoidal}.

From these assumptions, it follows that the Jacobian matrix of \eqref{eqn:SystemSFA}, given by $-\Delta+\pi A\frac{\partial \psi}{\partial x}(x)$, is Metzler. Therefore, the system \eqref{eqn:SystemSFA} is cooperative.

It is convenient to rewrite the system \eqref{eqn:SystemSFA} in the following form
\begin{equation}
\dot x =\Delta\left[-x+\pi H_1 \psi(x)\right],\quad x \in \R^n.
\label{eqn:SystemSFH}
\end{equation}
where the matrix $ H_1 := \Delta^{-1} A $. Denote also $H:=\pi H_1$.
Observe that it satisfies some useful properties:
\begin{itemize}
	\item $H$ is nonnegative and irreducible, so Theorem \ref{theorem:Perron-Frobenius} applies.
	\item All the row sums of $H$ are equal to $\pi$, that is $H \1_n = \pi \1_n$. It follows that $(\pi,\1_n)$ is the Perron-Frobenius eigenpair of $H$. 
	\item As the matrix $A$ is symmetrizable, also $H$ is symmetrizable, hence it has real eigenvalues.
\end{itemize}

We will see below that the existence (and the stability) of multiple equilibria is strictly related to the structure of the spectrum of the matrix $H$. 


\subsection{Existence of multiple equilibria: a necessary condition}
Consider the system \eqref{eqn:SystemSFH} (or \eqref{eqn:SystemSFA}) where each nonlinear function $\psi_i(x_i)$ satisfies the properties \eqref{assumption:1psiOdd}$\div$\eqref{assumption:4Sigmoidal}.
Let us start by recalling what is known for this system when we vary the parameter $ \pi$.
By construction, the origin is always an equilibrium point for \eqref{eqn:SystemSFH} (or \eqref{eqn:SystemSFA}).
When $\pi< 1$, $x=0$ is the only equilibrium point, it is globally asymptotically stable and locally exponentially stable.
This follows from diagonal dominance, and can be easily shown by a Lyapunov argument, see \cite{bib:FranciSrivastavaLeonard2015,UgoAbara2016Spectral,UgoAbara2015Existence}.
At $\pi=\pi_1= 1$ the system undergoes a pitchfork bifurcation, the origin becomes a saddle point and two more equilibria emerge, $x^+\in \R^n_+$ and $x^-\in \R^n_-$ \cite{bib:FranciSrivastavaLeonard2015}.
It follows from the analysis of \cite{UgoAbara2016Spectral} that, for all $\pi>1$, $x=0$ is an unstable equilibrium point, while both $x^+$ and $x^-$ are locally asymptotically stable with domain of attraction given by (at least) the entire orthant for any $ \pi >0$.
$\R^n_+$ (resp. $\R^n_-$) are in fact invariant for the system $\eqref{eqn:SystemSFH}$.
What happens outside these two orthants is however unknown.
When $\pi>1$ and $\pi-1$ sufficiently small, the behavior of the system \eqref{eqn:SystemSFH} outside $\R^n_+$ and $\R^n_-$ has been discussed in \cite{bib:FranciSrivastavaLeonard2015}. Only the three equilibrium points mentioned above are possible, two locally stable \cite{UgoAbara2016Spectral} and the origin as a saddle point.
However, when $\pi>1$ grows, the bifurcation analysis of \cite{bib:FranciSrivastavaLeonard2015} does not suffice anymore.

Our task is therefore to investigate the behavior of the system \eqref{eqn:SystemSFH} when $\pi>1$ grows and $x\in \R^n$ (case not described by \cite{UgoAbara2016Spectral} and \cite{bib:FranciSrivastavaLeonard2015}).
In particular, we would like to understand for what interval $ (1, \; \pi_2 )$ of the bifurcation parameter $ \pi $ extra equilibria not contained in $ \mathbb{R}^n_+/\mathbb{R}^n_-$ cannot appear, and what happens for $ \pi > \pi_2$.

The following theorem introduces a necessary condition that has to be verified in order to have an equilibrium point $\bar x$ in a generic orthant $\K\ne \R^n_+,\R^n_-$ for the system \eqref{eqn:SystemSFH}.

\begin{theorem}
\label{theorem:SFnecessaryCondition}	
Consider the system \eqref{eqn:SystemSFH} where each nonlinear function $\psi_i(x_i)$ satisfies the properties \eqref{assumption:1psiOdd}$\div$\eqref{assumption:4Sigmoidal}. 
If the system admits an equilibrium point $\bar x\in \K$, where $\K$ is an orthant in $\R^n$ and $\K \ne \R^n_+,\R^n_-$, then $\exists \,\lambda(H) \in \Lambda(H)$ such that $\lambda(H)>1$ and $\lambda(H) \ne \rho(H)$.
\end{theorem}
\proof
Let $\bar x \in \K$ be an equilibrium point for \eqref{eqn:SystemSFH}. Because $\Delta$ is diagonal and positive definite, from \eqref{eqn:SystemSFH} it follows
\begin{equation}
\bar x= H \psi(\bar x).
\label{eqn:eqPointK}
\end{equation}
First notice that if $\bar x \in \K$ also $\psi(\bar x)\in \K$, because $\psi_i(x_i)$ keeps the same sign of $x_i$ for all $i=1,2,\dots,n$.
Introduce the diagonal matrix $M(\bar x)=\diag\{m_1(\bar x_1),\dots,m_n(\bar x_n)\}$ where each element is given by
\begin{equation*}
m_i(\bar x_i)=\frac{\psi_i(\bar x_i)}{\bar x_i},\quad i=1,\dots,n.
\end{equation*} 
Since $\bar x>_{\K} 0$, the ratio is well-posed. The dependence of $M(\bar x)$ from $\bar x$ will be omitted from now on.
From \eqref{assumption:4Sigmoidal} one gets
$ \abs{\psi_i(x_i)}< \abs{x_i}\quad \forall i, \;x_i\ne 0 $, which leads to $m_i=\frac{\psi_i(\bar x_i)}{\bar x_i}\in (0,1)$ for all $i$.
Then $0<\diag (M)< \1_n$.
Applying the change of coordinates $\bar z=M^{1/2}\bar x$, from $\psi(\bar x)=M \bar x$, we get
\begin{equation}
\bar z=M^{1/2}H M^{1/2} \bar z.
\label{eqn:EigenvalueEqCongruentH}
\end{equation}
Eq.~\eqref{eqn:EigenvalueEqCongruentH} represents the eigenvalue equation for the matrix $M^{1/2}H M^{1/2}$, that is $(1,\bar z)$ is an eigenpair of $M^{1/2}H M^{1/2}$.
Furthermore $\bar z \in \K$ since, for each $i$, $\bar z_i=\sqrt{m_i}\,\bar x_i$ and $m_i$ is strictly positive.
Observe that $M^{1/2}H M^{1/2}$ is nonnegative and irreducible and let its eigenvalues be arranged in a nondecreasing order.
Theorem~\ref{theorem:Perron-Frobenius} states that $\lambda_n\left(M^{1/2}H M^{1/2}\right)=\rho\left(M^{1/2}H M^{1/2}\right)$ is real and positive and that its associated eigenvector is real and in the positive orthant of $\R^n$.
Then $\rho\left(M^{1/2}H M^{1/2}\right)> 1$, since the eigenvector associated to $1$ is $\bar z \in \K \ne \R^n_+, \R^n_-$.
To prove that $\exists\;\lambda(H)\in \Lambda(H)$ such that $\lambda(H)>1$ we proceed in steps.

\begin{itemize}
	\item Since the matrix $A$ is symmetrizable, let $A=D_A S_A$ where $D_A$ is a diagonal matrix with positive diagonal entries and $S_A$ is a symmetric matrix. Notice that the matrix $S_A$ is still irreducible, nonnegative and with null diagonal.		
	
	\item The matrix $H$, which can be written as $H=\pi\Delta^{-1}D_AS_A$, is symmetrizable. Define the matrix $D_H:=\pi\Delta^{-1}D_A$, diagonal with positive diagonal entries. Then $S_A=D_H^{-1}H$ is the symmetrization of $H$ while $D_H^{-1}$ is the symmetrizer of $H$. 
	
	\item Consider the matrix $\tilde H$ defined as
	\begin{equation*}
	\tilde H:=
	D_H^{-1/2}H D_H^{1/2}
	=D_H^{1/2}S_A D_H^{1/2}.
	\end{equation*}
	By construction it is symmetric, nonnegative, irreducible and similar to $H$.
	Because $H$ and $\tilde H$ have the same eigenvalues, it is just necessary to prove that $\exists\;\lambda(\tilde H)\in \Lambda(\tilde H)$ such that $\lambda(\tilde H)>1$ and $\lambda(\tilde H)\ne \rho(\tilde H)$.
	
	\item The matrices $M^{1/2}HM^{1/2}$ and $M^{1/2}\tilde H M^{1/2}$ are similar. Indeed
	\begin{align*}
	M^{1/2}H M^{1/2}
	&=D_H^{1/2}\left[M^{1/2}\tilde H M^{1/2}\right]D_H^{-1/2}
	\\&
	\sim M^{1/2} \tilde H M^{1/2}.
	\end{align*}
	Then they have the same eigenvalues and in particular, from equation \eqref{eqn:EigenvalueEqCongruentH}, it follows that $(1,D_H^{-1/2}\bar z)$ is an eigenpair of the matrix $M^{1/2}\tilde H M^{1/2}$. Since $D_H^{-1/2}\bar z\in \K$, it follows that $\exists\;k\ne n$ such that $\lambda_k(M^{1/2}\tilde H M^{1/2})=1$.
	
	\item The matrix $M^{1/2}\tilde H M^{1/2}$ is symmetric and $M^{1/2}$ is nonsingular so it is possible to apply Theorem~\ref{theorem:Ostrowski}. There exists a positive real number $\theta_k$ such that 
	\begin{gather*}
	\lambda_{k}(M^{1/2}\tilde H M^{1/2})
	=\theta_{k}\lambda_{k}(\tilde H)
	\end{gather*}
	and
	\begin{equation*}
	\lambda_1\left(M^{1/2}(M^{1/2})^T\right)
	\le \theta_k \le 
	\lambda_n\left(M^{1/2}(M^{1/2})^T\right)
	\end{equation*}
	where $\lambda_1(M)=m_{\min}$ and $\lambda_n(M)=m_{\max}$. Then $\theta_k \le m_{\max}<1$.
	Since $k\ne n$ is the index for which $\lambda_k(M^{1/2}\tilde HM^{1/2})=1$, it follows
	\begin{gather*}
	1=\lambda_k(M^{1/2}\tilde H M^{1/2})
	=\theta_k\lambda_k(\tilde H)<\lambda_k(\tilde H).
	\end{gather*}
	Because $k \ne n$, this implies the existence of an eigenvalue $\lambda(\tilde H)\in \Lambda(\tilde H)$ such that $\lambda(\tilde H)>1$ and $\lambda(\tilde H)\ne \rho(\tilde H)$.
	Consequently, since $H$ and $\tilde H$ are similar, there exists $\lambda(H)\in \Lambda(H)$ such that $\lambda(H)>1$ and $\lambda(H)\ne \rho(H)$.
\end{itemize}
\qed

When instead of \eqref{eqn:SystemSFH} the system \eqref{eqn:SystemSFA} is considered, then the results are less sharp, since they depend on the diagonal terms, which are not all identical as in \eqref{eqn:SystemSFH}.

\begin{corollary}
	\label{lemma:SFnecessaryCondition}	
	Consider the system \eqref{eqn:SystemSFA}, where each nonlinear function $\psi_i(x_i)$ satisfies the properties \eqref{assumption:1psiOdd}$\div$\eqref{assumption:4Sigmoidal}. If the system admits an equilibrium point $\bar x\in \K$, where $\K$ is an orthant in $\R^n$ and $\K \ne \R^n_+,\, \R^n_-$, then $\exists \,\lambda(A) \in \Lambda(A)$ such that $\lambda(A)>0$ and $\lambda(A)\ne\rho(A)$ for which $\pi\,\lambda(A)> \delta_{\min}$.
\end{corollary}
\proof
	Let the symmetrizable matrix $A$ be written as $A=DS_A$, where $D$ is a diagonal matrix with positive diagonal entries and $S_A$ is a symmetric matrix. Then $H=\pi\Delta^{-1}DS_A$.
	Define a new matrix
	\begin{equation*}
		\tilde H:=\Delta^{1/2}D^{-1/2}HD^{1/2}\Delta^{-1/2}=\pi\Delta^{-1/2}\tilde A\Delta^{-1/2}
	\end{equation*}
	where $\tilde A$ is defined as $\tilde A:=D^{-1/2}AD^{1/2}=D^{1/2}S_A D^{1/2}$.
	By construction, $\tilde H$ is symmetric, similar to $H$ and congruent to $\tilde A$, while $\tilde A$ is symmetric, similar to $A$ and congruent to $S_A$.
	Because $\tilde H$ and $\tilde A$ are both symmetric, it is possible to apply Theorem~\ref{theorem:Ostrowski}. To simplify the notation, let $\tilde S:=\sqrt{\pi}\Delta^{-1/2}$ and $\tilde H=\tilde S \tilde A \tilde S^T$. Therefore, there exists a positive real number $\theta_k$ such that the following conditions hold
	\begin{gather}
	\lambda_k(\tilde H)=\theta_k\lambda_k(\tilde A)
	\label{eqn:cond1LambdatildeA}
	\\
	\lambda_1(\tilde S \tilde S^T)\le \theta_k \le \lambda_n(\tilde S \tilde S^T)
	\label{eqn:cond2LambdatildeA}
	\end{gather}
	From Theorem~\ref{theorem:SFnecessaryCondition}, $\exists\; \bar k\ne n$ such that $\lambda_{\bar k}(H)>1$. It follows, by similarity, that $\lambda_k(\tilde H)>1$ where $k\ne n$.	
	Then the condition \eqref{eqn:cond1LambdatildeA} where $\theta_k>0$ and $k\ne n$, yields $\lambda_k(\tilde A)>0$ and $\lambda_k(\tilde A)\ne \rho(\tilde A)$.
	Moreover, since $\tilde S \tilde S^T=\pi\Delta^{-1}$, then
	\begin{equation*}
	\lambda_1(\tilde S \tilde S^T)=\frac{\pi}{\delta_{\max}}
	,\quad
	\lambda_n(\tilde S \tilde S^T)=\frac{\pi}{\delta_{\min}}.
	\end{equation*}
	From \eqref{eqn:cond2LambdatildeA}, \eqref{eqn:cond1LambdatildeA} and the result of Theorem~\ref{theorem:SFnecessaryCondition}, it follows that
	\begin{equation*}
		1< \frac{\pi}{\delta_{min}}\,\lambda_k(\tilde A)
	\end{equation*}
	i.e. $\pi\,\lambda_k(\tilde A)>\delta_{\min}$. 
	But $\tilde A$ and $A$ are similar, that is, they have the same eigenvalues. Then $\pi\,\lambda_k(A)>\delta_{\min}$, which concludes the proof.
\qed

It is possible to relax the assumption \eqref{assumption:4Sigmoidal}. Define for each $i$ the coefficients $\mu_i$ as
\begin{equation*}
\mu_i:=\max_{x_i}\left\{ \frac{\psi_i(x_i)}{x_i} \right\},\quad i=1,\dots,n
\end{equation*}
and then define
\begin{equation}
\mu:=\max_{i=1,\dots,n}\left\{\mu_i\right\}.
\label{eqn:maximumCoefficient}
\end{equation}
This means that the condition
\begin{equation*}
\abs{\psi_i(x_i)}\le  \mu \cdot \abs{x_i}
\end{equation*}
holds for each $i$ and $x_i\in \R$.

\begin{theorem}
	\label{theorem:SFnecessaryCondition_relaxCond}	
	Consider the system \eqref{eqn:SystemSFH} where each nonlinear function $\psi_i(x_i)$ satisfies the properties \eqref{assumption:1psiOdd}$\div$\eqref{assumption:3psiSaturated}. If the system \eqref{eqn:SystemSFH} admits an equilibrium point $\bar x>_\K 0$, where $\K$ is an orthant in $\R^n$ and $\K \ne \R^n_+$, $\K \ne \R^n_-$, then $\exists \,\lambda(H) \in \Lambda(H)$ such that $\lambda(H)\ge 1/\mu$ and $\lambda(H) \ne \rho(H)$, where $\mu$ is given by \eqref{eqn:maximumCoefficient}.
\end{theorem}
The proof of this theorem and of the following corollary are omitted as they are completely analogous to that of Theorem~\ref{theorem:SFnecessaryCondition} and Corollary~\ref{lemma:SFnecessaryCondition}.

\begin{corollary}
	\label{lemma:SFnecessaryCondition_relaxCond}
	Consider the system \eqref{eqn:SystemSFA} where each nonlinear function $\psi_i(x_i)$ satisfies the properties \eqref{assumption:1psiOdd}$\div$\eqref{assumption:3psiSaturated}. If the system \eqref{eqn:SystemSFA} admits an equilibrium point $\bar x\in \K$, where $\K$ is an orthant in $\R^n$ and $\K \ne \R^n_+,\, \R^n_-$, then $\exists \,\lambda(A) \in \Lambda(A)$ such that $\lambda(A)>0$ and $\lambda(A)\ne\rho(A)$, for which $\pi\,\lambda(A)\ge \frac{\delta_{\min}}{\mu}$, where $\mu$ is given by \eqref{eqn:maximumCoefficient}.
\end{corollary}

\begin{remark}
	Observe that if $\bar{x}\in\R^n$ is any equilibrium point for the system \eqref{eqn:SystemSFH} (or \eqref{eqn:SystemSFA}), also $-\bar{x}$ is an equilibrium point as well. Indeed, as by assumption \eqref{assumption:1psiOdd} $\psi(-\bar{x})=-\psi(\bar{x})$, it follows that $H\psi(-\bar{x})=-H\psi(\bar{x})=-\bar{x}$.
\end{remark}

\begin{remark}
The necessary conditions given by Theorem~\ref{theorem:SFnecessaryCondition} and Theorem~\ref{theorem:SFnecessaryCondition_relaxCond} imply that in order to have an equilibrium point $\bar x\in \mathcal K$ for the system \eqref{eqn:SystemSFA}, where $\K$ is an orthant in $\R^n$ and $\K \ne \R^n_+,\R^n_-$, the number of nodes in the network should be strictly greater than three. 
The next proposition in fact shows that it is impossible for $A$ (and hence for $H$) to have a second positive eigenvalue which differs from the Perron-Frobenius eigenvalue if $n \leq 3$, and is of independent interest.
\end{remark}

\begin{proposition}
	Let $n\le 3$ and $A\in\R^{n\times n}$ be a irreducible, symmetrizable, nonnegative matrix with null diagonal. Then $A$ cannot have two different real positive eigenvalues.
\end{proposition}
\proof
	The matrix $A$ is symmetrizable which implies that its eigenvalues are real. Let them be arranged in a nondecreasing order, that is $\lambda_n(A)>\lambda_{n-1}(A)\ge\dots\ge \lambda_1(A)$. It is always true that 
	\begin{equation*}
	\begin{system}
		\sum_{i=1}^n \lambda_i(A)={\rm Tr}(A)=0,\,\forall n\\
		\prod_{i=1}^n \lambda_i(A)=\det(A).
	\end{system}		 
	\end{equation*}
	Since $A$ is nonnegative and irreducible it is possible to apply Theorem~\ref{theorem:Perron-Frobenius} from which it follows that $\lambda_n(A)=\rho(A)>0$.
	Then
	\begin{itemize}
		\item If $n=2$, $A=\begin{bmatrix} 0 & a_{12}\\ a_{21} & 0 \end{bmatrix}$ and the conditions on its trace and determinant become
		\begin{equation*}
		\begin{system}
		\lambda_1(A)+\lambda_2(A)=0\\
		\lambda_1(A)\lambda_2(A)=-a_{12}a_{21}<0
		\end{system}	
		\end{equation*}
		which yields $\lambda_1(A)<0$.
		
		\item If $n=3$, %
		$A=
		\begin{bmatrix}
		0 & a_{12} & a_{13} \\
		a_{21} & 0 & a_{23} \\
		a_{31} & a_{32} & 0
		\end{bmatrix}$ %
		and the conditions on its trace and determinant become
		\begin{equation*}
		\begin{system}
		\lambda_1(A)+\lambda_2(A)+\lambda_3(A)=0\\
		\lambda_1(A)\lambda_2(A)\lambda_3(A)=a_{12}a_{23}a_{31}+a_{13}a_{21}a_{32}\ge 0
		\end{system}
		\end{equation*}
		which yields $\lambda_1(A),\,\lambda_2(A)<0$.
	\end{itemize}
	\qed


\subsection{A geometric necessary and sufficient condition}
\label{subsection:GeometricInterpretation}

The following lemma provides a geometric interpretation of the condition of Theorem~\ref{theorem:SFnecessaryCondition}. 
Consider the Laplacian $ L_1 = I - H_1 $. 
Since $ \rho(H_1) =1$, by construction, the least eigenvalue $ \lambda_1 (L_1)= 1 - \rho(H_1) $ is the origin. 
Recall that the second leftmost eigenvalue of $L_1 $, $ \lambda_{2} (L_1) $, is called the algebraic connectivity of $ L_1 $ \cite{Olfati2003Consensus}. 
If $ \lambda_{n-1}(H_1) $ is the second largest eigenvalue of $ H_1 $, then the algebraic connectivity of $L_1 $ is  $ \lambda_2 (L_1) = 1- \lambda_{n-1} (H_1) $.

\begin{lemma}
\label{lemma:alg_conn}
The range of values of $ \pi $ for which no extra equilibrium of the system \eqref{eqn:SystemSFH} (other than $0$, $ x^+ $ and $ x^-$) can appear is given by $ (1, \; \pi_2) $, with $ \pi_2 = \frac{1}{\lambda_{n-1}(H_1)}>1$, and hence it is determined by the algebraic connectivity $ \lambda_2 (L_1)$.
\end{lemma}

\proof
For \eqref{eq:PF_equal1} it holds that $ \sum_j h_{1, ij} = \rho(H_1) = 1 $ $ \forall \; i $.
This means that the Laplacian $ L_1  $ has all identical Ger\v{s}gorin disks, 
all centered at 1 and passing through the origin:
\[
\left\{ s \in \mathbb{C} \text{  s.t.  } | s - 1| \leqslant \sum_{j=1}^n  \frac{a_{ij}}{\delta_i } =1 \right\}.
\]
From Ger\v{s}gorin Theorem \cite{bib:HornJohnson1990}, the eigenvalues of $L_1$ are located in the union of the $n$ disks. 
The least eigenvalue $ \lambda_1 (L_1)=0$ lies on the boundary of the Ger\v{s}gorin disks.
The other eigenvalues of $ L_1 $ are strictly inside the disks because of irreducibility and positive semidefiniteness of $L_1$.
When $ \pi > 1 $, all the eigenvalues of $H = \pi H_1 $ are increased in modulus, and the Ger\v{s}gorin disks of $I - H $ are still all centered in 1 but have radius $ \pi $.
The condition of Theorem~\ref{theorem:SFnecessaryCondition}, $ \exists \, \lambda(H) >1$, $ \lambda(H) \neq \rho(H) $ corresponds to the existence of $  \pi $ s.t. $ \pi \lambda(H_1) > 1 $ for $ \lambda (H_1 ) > 0$ and $ \lambda (H_1 ) \neq 1$,  i.e., it corresponds to requiring that $ \pi > \pi_2 $, where $ \pi_2 = \frac{1}{\lambda_{n-1} (H_1)}>1 $. 
When this happens, $  1 - \pi \lambda_{n-1} (H_1 ) $ crosses into the left half of the complex plane.  
\qed

The insight given by Lemma~\ref{lemma:alg_conn} allows to show that the condition of Theorem~\ref{theorem:SFnecessaryCondition} is also sufficient. Since we use bifurcation theory and singularity analysis in the proof, we have however to restrict to the case of algebraic connectivity $ \lambda_2 (L_1) $ which is a simple eigenvalue (property which is generic for weighted graphs, see \cite{Poignard2017Spectra}).

\begin{theorem}
\label{theorem:SFnecessaryandSufficientCondition}	
Consider the system \eqref{eqn:SystemSFH} where each nonlinear function $\psi_i(x_i)$ satisfies the properties \eqref{assumption:1psiOdd}$\div$\eqref{assumption:4Sigmoidal}. 
Assume further that the second largest eigenvalue of $ H_1 $, $ \lambda_{n-1}(H_1) $ is simple.
The system admits an equilibrium point $\bar x\in \K$, where $\K$ is an orthant in $\R^n$ and $\K \ne \R^n_+,\R^n_-$, if and only if $ \pi > \pi_2 = \frac{1}{\lambda_{n-1}(H_1)}>1$.
\end{theorem}
\proof
Necessity was proven in Theorem~\ref{theorem:SFnecessaryCondition}. In fact, as shown in Lemma~\ref{lemma:alg_conn}, since $ \rho(H_1) =1$, the condition of Theorem~\ref{theorem:SFnecessaryCondition} corresponds to $ \pi  > \pi_2 = \frac{1}{\lambda_{n-1} (H_1)}>1$.
To show that mixed sign equilibria appear exactly at $ \pi = \pi_2 = \frac{1}{\lambda_{n-1} (H_1)} $, we use bifurcation theory, in particular singularity theory and Lyapunov-Schmidt reduction for pitchfork bifurcations \cite{golubitsky2000singularities}, see also \cite{bib:FranciSrivastavaLeonard2015}. 
Consider 
\beq
\Phi(x, \, \pi) = - x + \pi H_1 \psi(x) =0 .
\label{eq:signularity_1}
\eeq
Denote $ v_2 $ and $w_2$, $ w_2^T v_2 =1 $, the right and left Fiedler vectors, i.e., the eigenvectors relative to $ \lambda_{n-1} (H_1 ) $. 
From Theorem~\ref{theorem:Perron-Frobenius}, $ v_2,w_2 \in \K $ for some $ \K \ne \R^n_+,\R^n_-$.
Let $J=\pde{\Phi}{x}(0,\,\pi_2)=-I+\pi_2 H_1$. Observe that $\text{range}(J)=\left(\ker (J ^T)\right)^\perp =\left({\rm span}\{ w_2\}\right)^\perp$ and that $\ker(J)={\rm span}\{ v_2\}$.
Let $E$ denote the projection of $\R^n$ onto ${\rm range}(J)=\left({\rm span}\{ w_2\}\right)^\perp$, $E=I-v_2w_2^T$, and $I-E$ the projection onto $\ker(J)={\rm span}\{ v_2\}$.
Split $ x $ accordingly: $ x = y v_2 + r $, where $ y \in \mathbb{R} $ and $ r = E x \in ( {\rm span} \{ w_2 \})^\perp$.

Near $ (0, \, \pi_2 ) $, \eqref{eq:signularity_1} can be split into 
\begin{subequations}
\begin{align}
E \Phi(x, \, \pi) & = E (- x + \pi H_1 \psi(x) )=0 \label{eq:signularity_2a}\\
(I-E) \Phi(x, \, \pi) & = (I-E) (- x + \pi H_1 \psi(x) )=0 \label{eq:signularity_2b}
\end{align}
\end{subequations}
Since $ \lambda_{n-1}(H_1) $ is simple, at $ \pi= \pi_2 $ \eqref{eq:signularity_1} has a simple singularity. Hence for \eqref{eq:signularity_2a} the implicit function theorem applies and it is possible to express $ r = R(y v_2, \, \pi) $. The explicit expression of $ R(\cdot ) $ is not needed in what follows. 
Replacing $ r$ in \eqref{eq:signularity_2b}, we get the center manifold
\beq
\begin{split}
 \phi(y, \, \pi) = & (I-E) \big(-y v_2 - R(y v_2 , \, \pi) \\
 &  + \pi H_1 \psi(y v_2 + R(y v_2, \, \pi ) )\big )=0 .
\end{split}
\label{eq:signularity_3}
\eeq
Define $g(y,\pi)=w_2^T\phi(y,\pi)$, where $w_2\in \left(\text{range}\,J\right)^\perp$. The recognition problem for a pitchfork bifurcation requires to compute at $ (0, \, \pi_2 ) $ the partial derivatives $ g_y $, $ g_{yy} $, $ g_{yyy} $, $ g_\pi $, $ g_{\pi y} $.
In this case the calculation is simplified by the fact that $ \Phi(x, \, \pi) $ is odd in $ x$.  
For instance, since $ \Phi $ has a singularity at $ (0, \, \pi_2) $ it follows that the directional derivative along $ v_2 $ vanishes:
\[
\begin{split}
\Phi_y(0, \, \pi_2) & = \left. \pde{\Phi}{x} (x, \, \pi) \right|_{(0, \, \pi_2)}   \! \! \! \! \! \! \! \! \!  v_2 = \left. \left(-I + \pi H_1 \pde{\psi(x)}{x}\right)  \right|_{(0, \, \pi_2)} \! \! \! \! \! \! \! \! \! \! \!  v_2 \\
& = (-I + \pi_2 H_1) v_2 = -v_2 + \pi_2 \lambda_{n-1}(H_1) v_2 =0
\end{split}
\]
where we have used $ \left. \pde{\psi(x_i)}{x_i} \right|_{0} =1 $. 
Similarly (see \cite{golubitsky2000singularities})
\[
\begin{split}
\Phi_{yy} (0, \, \pi_2) & =  \pde{}{x}\left( \left. \pde{\Phi}{x} (x, \, \pi)\right|_{(0, \, \pi_2)}   \! \! \! \! \! \! \! v_2 \right) v_2 \\
&
=\pi H_1 \pde{}{x}\left(  \left. \pde{\psi(x)}{x} v_2 \right) \right|_{ (0, \, \pi_2)} \! \! \! \! \! \! \! \! \! \! v_2 =0
\end{split}
\]
because $ \left. \ppde{\psi(x_i)}{x_i} \right|_{0} =0 $, and 
\[
\Phi_\pi(0, \, \pi_2)  = \left. \pde{\Phi}{\pi} (x, \, \pi) \right|_{(0, \, \pi_2)}  \! \! \! \! \! \! = \left.  H_1 \psi(x)  \right|_{(0, \, \pi_2)} =0
\]
since $ \psi_i(0)=0$.
The two remaining partial derivatives are
\[
\begin{split}
\Phi_{\pi y} (0, \, \pi_2) & =  \pde{}{\pi}\left( \left. \pde{\Phi}{x} (x, \, \pi)\right|_{(0, \, \pi_2)}  v_2 \right)  \\
& = H_1  \left. \pde{\psi(x)}{x}  \right|_{ (0, \, \pi_2)}\! \! \! \! \! \!  v_2 = H_1 v_2 = \lambda_{n-1}(H_1) v_2 
\end{split}
\]
and (using a notation similar to \cite{golubitsky2000singularities}, (3.16))
\[
\begin{split}
\Phi_{yyy} (0, \, \pi_2) & =   \left. \frac{\partial ^3 \Phi}{\partial x^3} (x, \, \pi)\right|_{(0, \, \pi_2)} ( v_2 , \, v_2 , \, v_2 ) \\
& =\pi_2 H_1 \begin{bmatrix} \beta_1 & \\ & \ddots \\ & & \beta_n \end{bmatrix} 
\begin{bmatrix} v_{2,1}^3 \\ \vdots \\ v_{2,n}^3 \end{bmatrix}
\end{split}
\]
where, from (A.4), $ \beta_i = \left. \frac{\partial ^3 \psi_i(x_i)}{\partial x_i^3} \right|_{0} <0$. 
Consequently, for the projections along $ v_2$: $ \phi_y(0, \, \pi_2)= (I-E) \Phi_y(0, \, \pi_2)=0 $, and, similarly, $ \phi_{yy}(0, \, \pi_2)= 0 $, $ \phi_\pi(0, \, \pi_2)=0 $.
As for the two nonvanishing partial derivatives: 
$
\phi_{\pi y}(0, \, \pi_2) = (I-E) \Phi_{\pi y} (0, \, \pi_2) = \lambda_{n-1}(H_1)v_2
$ 
and
$
\phi_{yyy}(0, \, \pi_2) = (I-E) \Phi_{yyy} (0, \, \pi_2) = v_2\sum_{i=1}^n \beta_i w_{2,i} v_{2,i}^3
$. 
Therefore $ g_y(0, \, \pi_2)= g_{yy}(0, \, \pi_2)= g_\pi(0, \, \pi_2)=0 $ and $g_{\pi y}(0, \, \pi_2) = \lambda_{n-1}(H_1)>0$, $g_{yyy}(0, \, \pi_2) = \sum_{i=1}^n \beta_i w_{2,i} v_{2,i}^3<0$ since $w_{2,i}v_{2,i}^3\ge 0$ $\forall i$, which completes the recognition problem for a pitchfork bifurcation \cite{golubitsky2000singularities}.
Hence at $ \pi = \pi_2 $ the system crosses a second bifurcation through the origin and two new equilibria appear along $ {\rm span}\{v_2\}$.
Since $ v_2 \in \mathcal{K} $, these equilibria must belong to $ \mathcal{K} $ and $ - \mathcal{K}$. 
\qed

Notice that unlike most arguments based on singularity analysis of bifurcations, our result in Theorem~\ref{theorem:SFnecessaryandSufficientCondition} is not a local one, as our proof of necessity (Theorem~\ref{theorem:SFnecessaryCondition}) is nonlocal. 

When $ \lambda_{n-1}(H_1) $ has multiplicity higher than one, then singularity analysis based on the normal form of a pitchfork bifurcation does not apply, although we expect that similar sufficiency results can be obtained through more advanced bifurcation theory. 

A pictorial view of the situation described in Theorem~\ref{theorem:SFnecessaryandSufficientCondition} will be given below in Example~\ref{ex:example1}, see in particular Fig.~\ref{fig:example_n6} and Fig.~\ref{fig:example_n6disks} (a), (b).

If instead we look at system \eqref{eqn:SystemSFA}, and at the Laplacian $ L= \Delta - A$, then, when $\pi = 1 $, the Ger\v{s}gorin disks are centered at $ \delta_i $ and have different radii, equal to $ \delta_i $. 
However, this cannot be straightforwardly reformulated in terms of $ \rho(A) $, as \eqref{eq:PF_disequal1} (instead of \eqref{eq:PF_equal1}) now holds: $ \delta_{\min} \leq \rho(A) \leq \delta_{\max}$. 
When exploring the values $ \pi >1$, then the Ger\v{s}gorin disks of $\tilde L = \Delta-\pi A$ are contained one in the other, according to the corresponding $\delta_i$, and all have nonempty intersection with the left half of $\C$, see Fig.~\ref{fig:example_n6disks}, panels (c) and (d).
In this case, while the condition introduced by Theorem~\ref{theorem:SFnecessaryCondition} represents the necessary and sufficient condition for the negativity of a second leftmost eigenvalue of $\tilde L$, the condition introduced by Corollary~\ref{lemma:SFnecessaryCondition} is just necessary, due to the presence of the diagonal matrix $\Delta$. Using a similar reasoning it is in fact possible to prove that $\lambda_{n-1}(H)=\pi \theta \lambda_{n-1}(A)$, where the value of the positive constant $\theta$ is not fixed but $\delta_{\min}\le \theta\le \delta_{\max}$.
	


\subsection{Stability properties of multiple equilibria}

\begin{theorem}
	\label{theorem:SFstability}
	Suppose the system \eqref{eqn:SystemSFH} admits an equilibrium point $\bar x\in \K$, where $\K$ is an orthant in $\R^n$ and $\K \ne \R^n_+,\,\R^n_-$.
	If
	\begin{equation}
	\pi\,\max_j\left\{\frac{\partial \psi_j}{\partial x_j}(\bar x_j)\right\}<1
	\label{eqn:ASstability}
	\end{equation}
	then $\bar x$ is locally asymptotically stable.
	Instead, if 
	\begin{equation}
	\pi\min_j\left\{\frac{\partial \psi_j}{\partial x_j}(\bar x_j)\right\}>1
	\label{eqn:unstability}
	\end{equation}
	then the equilibrium point $\bar x$ is unstable.
\end{theorem}
\proof
	Let $\bar x \in \K$ be an equilibrium point for \eqref{eqn:SystemSFH}.
	To study the behavior of the dynamical system \eqref{eqn:SystemSFH} near the equilibrium point $\bar x$, consider the linearization around $\bar x$
	\begin{equation}
	\dot x=\Delta\left[-I+H\;\frac{\partial \psi}{\partial x}(\bar x)\right](x-\bar x).
	\label{eqn:Linearization}
	\end{equation} 
	Under the condition \eqref{eqn:ASstability} it can be proven that the matrix $I-H\;\frac{\partial \psi}{\partial x}(\bar x)$ is strictly diagonally dominant, that is $1>\sum_{j\ne i} h_{ij}\,\frac{\partial \psi_j}{\partial x_j}(\bar x_j)$ for all $i=1,\dots,n$. Indeed
	\begin{align*}
		\sum_{j\ne i} h_{ij}\,\frac{\partial \psi_j}{\partial x_j}(\bar x_j)
		&\le \max_j\left\{\frac{\partial \psi_j}{\partial x_j}(\bar x_j) \right\}\,\sum_{j\ne i} h_{ij}
		\\&= \pi \,\max_j\left\{\frac{\partial \psi_j}{\partial x_j}(\bar x_j)\right\} < 1.
	\end{align*}
	Then it is possible to apply Theorem~\ref{theorem:SDDmatrix} and state that, under the condition \eqref{eqn:ASstability}, each eigenvalue of the matrix $I-H\;\frac{\partial \psi}{\partial x}(\bar x)$ has strictly positive real part.
	Therefore, each eigenvalue of the matrix $\Delta\left(-I+H\;\frac{\partial \psi}{\partial x}(\bar x)\right)$ has strictly negative real part, i.e. $\bar x$ is locally asymptotically stable.

	Now suppose that \eqref{eqn:unstability} holds and consider the linearization around $\bar x$ \eqref{eqn:Linearization}.
	Let $\bar H:=H\;\frac{\partial \psi}{\partial x}(\bar x)$ and notice that it is nonnegative and irreducible (therefore it is possible to apply Theorem~\ref{theorem:Perron-Frobenius}). 
	If the matrix $\bar H-I$ admits an eigenvalue with positive real part, it is possible to conclude that the equilibrium point $\bar x$ is unstable. 
	For each $i=1,\dots,n$, it holds that
	\begin{align*}
	\sum_{j=1}^n \bar h_{ij}
	& =\sum_{j=1}^n h_{ij}\,\frac{\partial \psi_j}{\partial x_j}(\bar x_j) \ge \min_j\left\{\frac{\partial \psi_j}{\partial x_j}(\bar x_j)\right\}\, \sum_{j=1}^n h_{ij}\\
	&= \min_j\left\{\frac{\partial \psi_j}{\partial x_j}(\bar x_j)\right\}\, \pi > 1
	\end{align*}
	under the hypothesis \eqref{eqn:unstability}.
	According to Corollary~\ref{corollary:Perron-Frobenius}, $\rho(\bar H)\ge \min_i \left(\sum_{j=1}^n \bar h_{ij}\right)$. From the previous reasoning, it follows that $\rho(\bar H)>1$.
	Therefore the matrix $\bar H-I$ admits a real positive eigenvalue given by $\rho(\bar H)-1$, which implies that the equilibrium point $\bar x$ is unstable.
\qed


\section{Location of the mixed-sign equilibria}
\label{section:EquilibriumPointPositiveOrthant}
In this Section we restrict our analysis to the special case of all identical $\psi_i(x_i)$. In this case the equilibrium point in the positive orthant has all identical components as shown in the following lemma.
\begin{lemma}
	\label{lemma:NormEqPos}
	Consider the system \eqref{eqn:SystemSFH} where each $\psi_i(x_i)$ satisfies the properties \eqref{assumption:1psiOdd}$\div$\eqref{assumption:4Sigmoidal} and $\psi_i(\xi)=\psi_j(\xi)$ $\forall \, i,j$ and $\forall \, \xi \in \R$.
	When $\pi>1$ the positive equilibrium point $x^+\in\R^n_+$ is such that $\frac{\psi_i(x_i^+)}{x_i^+}=\frac{1}{\pi}$.
	Furthermore $x^+$ is the Perron-Frobenius (right) eigenvector of $H$. 
\end{lemma}
\proof
	For the matrix $H$ we have that the Perron-Frobenius eigenvalue is $\rho(H)=\pi>1$ and the corresponding eigenvector is $\1_n$.
	At $x^+$ it must be $H\psi(x^+)=x^+$. If $\psi(x^+)=\frac{1}{\pi}\,x^+$ with $x^+=\alpha\1_n$ for some $\alpha\in (0,1)$ then we have $H\psi(x^+)=\frac{1}{\pi}\,\alpha H\1_n=\alpha\1_n$ from which we get the eigenvalue/eigenvector equation $H\1_n=\pi\1_n$. It follows that $x^+=\alpha\1_n$ is the Perron-Frobenius eigenvector of $H$.
	The specific value of $\alpha$ depends on the functional form of $\psi_i(x_i)$.
\qed

A more important consequence is that for each value of $\pi$ the positive equilibrium point $x^+$ provides an upper bound on the norm that any mixed-sign equilibrium $ \bar x $ can assume.

\begin{theorem}
	\label{theorem:NormBound} 
	Consider the system \eqref{eqn:SystemSFH}, where each nonlinear function $\psi_i(x_i)$ satisfies the properties \eqref{assumption:1psiOdd}$\div$\eqref{assumption:4Sigmoidal} and $\psi_i(\xi)=\psi_j(\xi)$ $\forall\,  i,j$ and $\forall \, \xi \in \R$. If the system admits an equilibrium point $\bar x\in \K$, where $\K$ is an orthant in $\R^n$ and $\K \ne \R^n_+,\,\R^n_-$, then $\norm{\bar x}\le \norm{x^+}$.
\end{theorem}
In order to prove the theorem we use the following lemma.
\begin{lemma}
	\label{lemma:NormEqbar}
	Consider the system \eqref{eqn:SystemSFH}, where each nonlinear function $\psi_i(x_i)$ satisfies the properties \eqref{assumption:1psiOdd}$\div$\eqref{assumption:4Sigmoidal} and $\psi_i(\xi)=\psi_j(\xi)$ $\forall \, i,j$ and $\forall \, \xi \in \R$. If the system admits an equilibrium point $\bar x\in \K$, where $\K$ is an orthant in $\R^n$ and $\K \ne \R^n_+,\, \R^n_-$, then $\frac{\psi_i(\bar x_i)}{\bar x_i}\ge \frac{1}{\pi}$ for all $i$.
\end{lemma}
\proof
	Introduce the diagonal matrix $M(\bar x)=\diag\{m_1(\bar x_1),\dots,m_n(\bar x_n)\}$ as in Theorem~\ref{theorem:SFnecessaryCondition}, where each element is given by
	\begin{equation*}
		m_i(\bar x_i)=\frac{\psi_i(\bar x_i)}{\bar x_i},\quad i=1,\dots,n
	\end{equation*} 
	Since $\bar x>_{\K} 0$ the ratio is well-posed. The dependence of $M(\bar x)$ from $\bar x$ will be omitted from now on. 
	Let $m_{\min}=\min_i\left\{\frac{\psi_i(\bar x_i)}{\bar x_i}\right\}$ for all $i$ and $l\in\{1,\dots,n\}$ be the index such that $\frac{\psi_l(\bar x_l)}{\bar x_l}=m_{\min}$ and suppose without loss of generality that $\bar x_l>0$. By definition of $m_{\min}$, $\psi_i(\bar x_i)\le \psi_l(\bar x_l)$ and $\bar x_i\le \bar x_l$ for all $i\ne l$.
	At $\bar x$ it must be $\bar x=H\psi(\bar x)$, that is $\bar x_i=\sum_{j\ne i} h_{ij}\psi_j(\bar x_j)$ for all $i=1,\dots,n$.
	This yields
	\begin{equation*}
		\bar x_l
		=\sum_{j\ne l} h_{lj}\psi_j(\bar x_j)
		\le \left(\sum_{j\ne l} h_{lj}\right)\psi_l(\bar x_l)
		=\pi\, \psi_l(\bar x_l)
	\end{equation*}
	which implies that $m_{\min}=\frac{\psi_l(\bar x_l)}{\bar x_l}\ge \frac{1}{\pi}$.
	Then, since $m_i\ge m_{\min}$ for all $i\ne l$, it follows that $\frac{\psi_i(\bar x_i)}{\bar x_i}\ge\frac{1}{\pi}$ for all $i=1,\dots, n$.
\qed

\proof[Proof of Theorem~\ref{theorem:NormBound}]
	Let $x^+\in \R^n_+$ and $\bar x\in \K$ be equilibrium points for the system \eqref{eqn:SystemSFH}.
	From Lemma~\ref{lemma:NormEqPos} and Lemma~\ref{lemma:NormEqbar} it follows that $\frac{\psi_i(\bar x_i)}{\bar x_i}\ge \frac{1}{\pi}=\frac{\psi_i(x_i^+)}{x_i^+}$ for all $i=1,\dots,n$. This implies that $\abs{\bar x_i}\le x_i^+$ for each $i$, which yields $\abs{\bar x}\le x^+$. Therefore $\norm{\bar x}\le \norm{x^+}$.
\qed


\section{Numerical Analysis}
\label{section:NumericalAnalysis}
In this Section we first look at the trajectories of a specific numerical example of size $ n=6 $. 
Then we perform a computational analysis of the properties of the equilibria for a system of size $ n = 20 $.

\begin{example}
\label{ex:example1}
Consider a network of $n=6$ nodes of weighted adjacency matrix
\begin{equation*}
A=\left[
\arraycolsep=2.0pt
\begin{array}{cccccc}
0 &0.1477 &0 &0.1698 &0 &0.0135 \\
0.4242 &0 &0.2626 &0.3621 &0 &0 \\
0 &0.1889 &0 &0 &0.2502 &0.4158 \\
0.4036 &0.2997 &0 &0 &0 &0 \\
0 &0 &0.2427 &0 &0 &0.2513 \\
0.0301 &0 &0.4474 &0 &0.2787 &0
\end{array}
\right]
\end{equation*}
which is irreducible, nonnegative and symmetrizable.
Consider the system \eqref{eqn:SystemSFA} (and \eqref{eqn:SystemSFH}) and assume that each nonlinear function $\psi_i(x_i)$ is given by the Boltzmann function $\psi_i(x_i)=\frac{1-e^{-2x_i}}{1+e^{-2x_i}}$, which satisfies the assumptions \eqref{assumption:1psiOdd}$\div$\eqref{assumption:4Sigmoidal}.
\begin{figure*}[ht!]
\centering
\subfloat[]{
\includegraphics[trim=0cm 0cm 0cm 0cm, clip=true,  width=5cm]{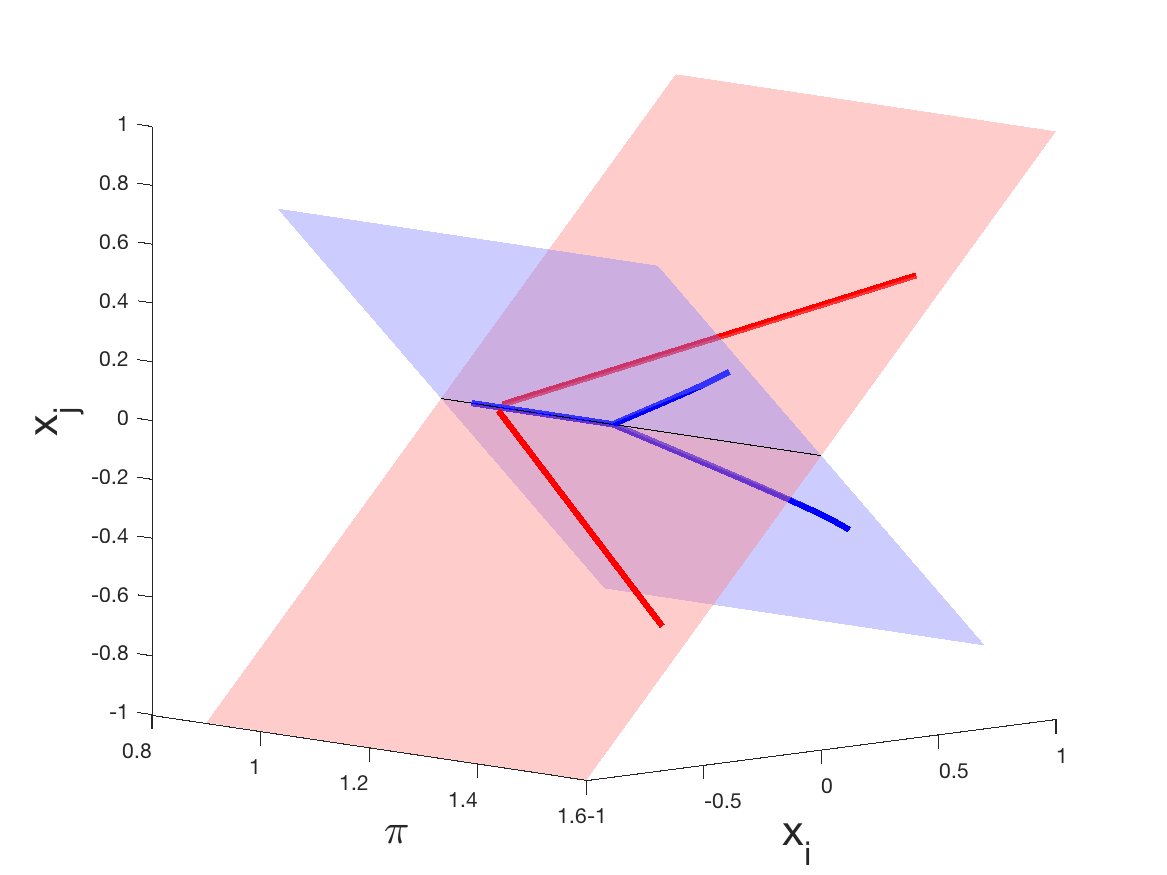}}
\subfloat[]{
\includegraphics[trim=0cm 0cm 0cm 0cm, clip=true,  width=5cm]{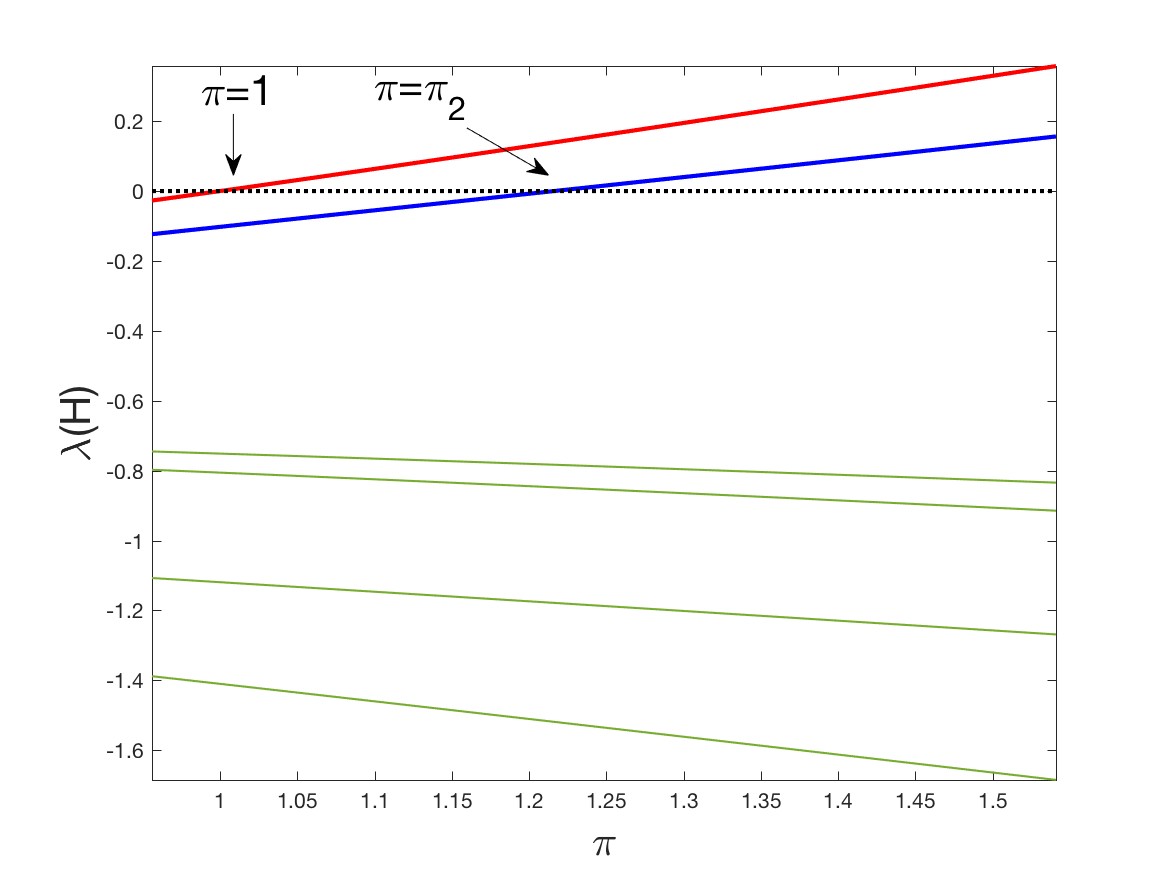}}
\subfloat[]{
\includegraphics[trim=0cm 0cm 0cm 0cm, clip=true,  width=5cm]{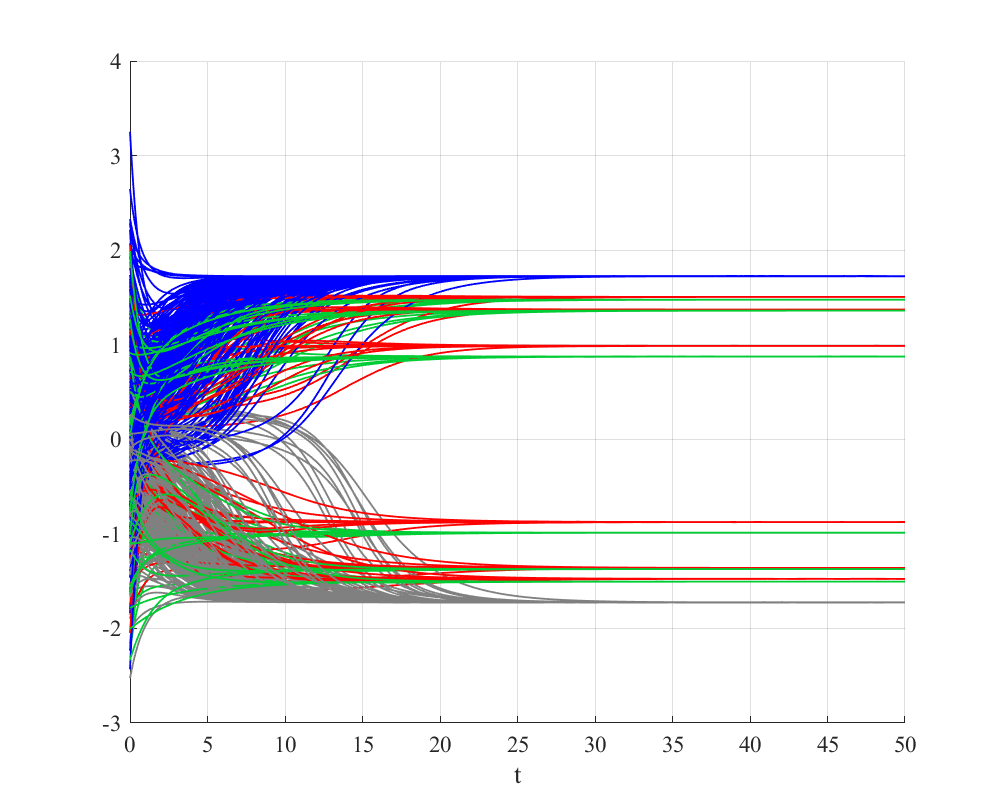}}
\caption{\small Example~\ref{ex:example1}. (a): Bifurcation diagram near $ x=0$ for two components $ x_j $ and $ x_j $. When $ \pi =1$ then the origin bifurcates a first time, and the two locally stable equilibria  $ x^{\pm} $ (red curve) are created along the consensus manifold (red plane). When $ \pi = \pi_2 = \frac{1}{\lambda_{n-1} (H_1)}$ a second bifurcation occurs in the origin, with 2 equilibria (blue curve) locally branching out along $ {\rm span}\{ v_2 \}$ (blue plane). (b): The corresponding set of eigenvalues of $ -I + \pi H_1 $ (linearization in the origin). The two crossing at $ \pi =1 $ and $ \pi= \pi_2 $ are highlighted.  (c): $100$ random initial conditions converging to $x^+$ (blue), $x^-$ (gray), $\bar x_3$ (red) and $-\bar x_3$ (green).}
\label{fig:example_n6}
\end{figure*}
The Perron-Frobenius eigenvalue of the matrix $A$ is $\rho(A)=0.706$ and its second largest positive eigenvalue is $\lambda_{n-1}(A)=0.515$, meaning $ \lambda_{n-1} (H_1) = 0.822$ or $ \pi_2 =1.216$.
Hence we have:
\begin{itemize}
\item for $ \pi<1 $, $ \bar x =0 $ is the only equilibrium point;
\item for $ \pi \in (1, \, 1.216) $ the only equilibria are $ \{ 0, \, x^+, \, x^- \}$, with $ x^\pm = \pm \alpha \1$, see red curve in Fig.~\ref{fig:example_n6} (a);
\item for $ \pi > 1.216 $ unstable equilibria $ \bar x_2 \in \K=\diag\{-1,\,-1,\,1,\,-1,\,1,\,1\} $ and $- \bar x_2 \in - \K $ bifurcate from $0$, see blue curve Fig.~\ref{fig:example_n6} (a).
\end{itemize}
When $ \pi $ grows further, new bifurcation points soon appear. 
Fig.~\ref{fig:example_n6} (b) shows that these bifurcations are not associated to singularities in the origin, but rather they branch out of $ \pm \bar x_2 $.
For instance choosing $ \pi=1.838$, also the condition $\pi\lambda_{n-1}(A)>\delta_{\min}$ is satisfied.
Numerical computations confirm the existence of at least 3 equilibria in $ \K$, denoted $\bar x_3$, $\bar x_4$, $\bar x_5$,
which shows that more than one equilibrium can exist per orthant.
The condition presented in Theorem~\ref{theorem:SFstability} is satisfied for $\bar x_3$, and this assures its local asymptotic stability.
Instead, $\bar x_4$ and $\bar x_5$ are unstable.
This is shown in the simulation of Fig.~\ref{fig:example_n6} (c).
Following the reasoning introduced in Section \ref{subsection:GeometricInterpretation}, it is possible to observe that at $ \pi= 1.838$ the second leftmost eigenvalue of the matrix $\tilde L=\Delta-\pi A$, $\lambda_2(\tilde L)=-0.302$, is negative. However, notice that even if the necessary condition presented in Corollary~\ref{lemma:SFnecessaryCondition} is satisfied, the condition $\pi \lambda_{n-1}(A)>\delta_{\max}$ does not hold. 
Fig.~\ref{fig:example_n6disks} shows the Ger\v{s}gorin's disks and the eigenvalues of different matrices, respectively $I-H_1$ (when $\pi=1$), $I-\pi H_1$ (when $\pi>1$), $L=\Delta-A$ (when $\pi=1$), $\tilde L=\Delta-\pi A$ (when $\pi>1$). Notice that while it is possible to determine the exact value of $\lambda_2(I-\pi H_1) = 1 - \pi \lambda_{n-1} (H_1)$, it is only possible to give a bound for the second leftmost eigenvalue of $\tilde L$, i.e. $\delta_{\max}[1-\lambda_{n-1}(H)] \le \lambda_2(\tilde L)\le \delta_{\min} [1-\lambda_{n-1}(H)]$ or $\delta_{\min}-\pi\lambda_{n-1}(A) \le \lambda_2(\tilde L)\le \delta_{\max}-\pi\lambda_{n-1}(A)$. Finally, observe that this eigenvalue has to be in the union of all the Ger\v{s}gorin's disks of $\tilde L$, but it is not possible to define a priori the smallest disk that contains $\lambda_2(\tilde L)$.
\begin{figure*}[ht!]
	\centering
	\subfloat[]{
		\includegraphics[trim=0cm 0cm 0cm 0cm, clip=true,  width=4.1cm]{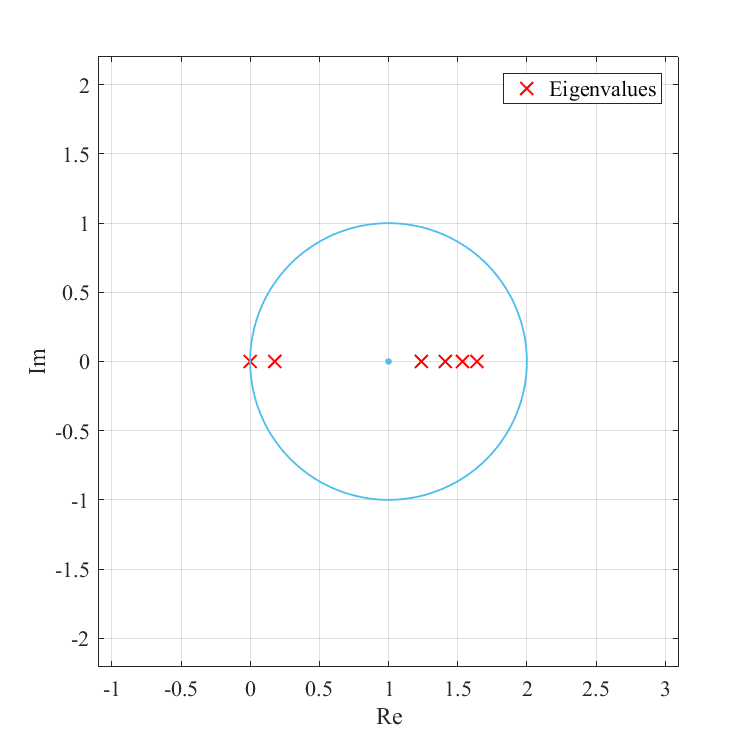}}
	\subfloat[]{
		\includegraphics[trim=0cm 0cm 0cm 0cm, clip=true,  width=4.1cm]{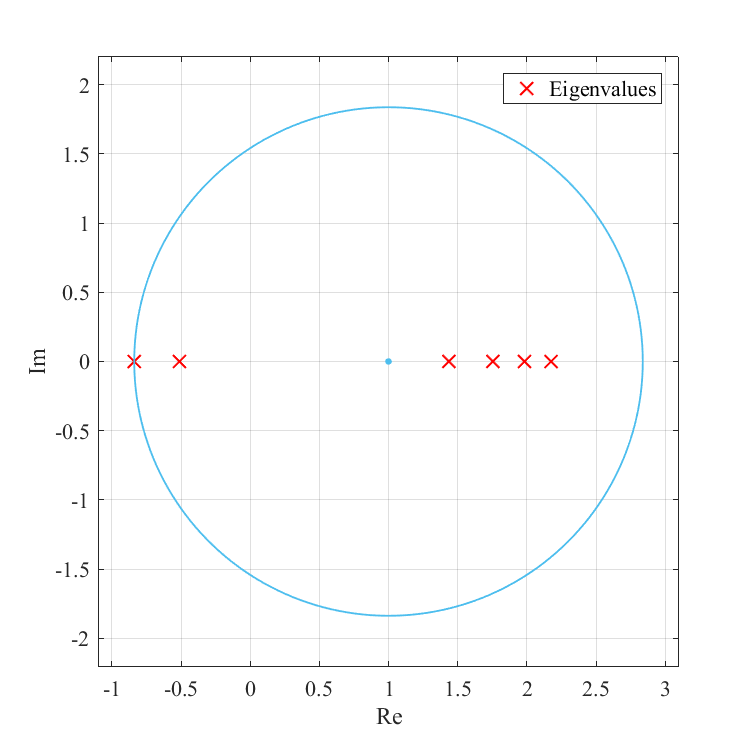}}
	\subfloat[]{
		\includegraphics[trim=0cm 0cm 0cm 0cm, clip=true,  width=4.1cm]{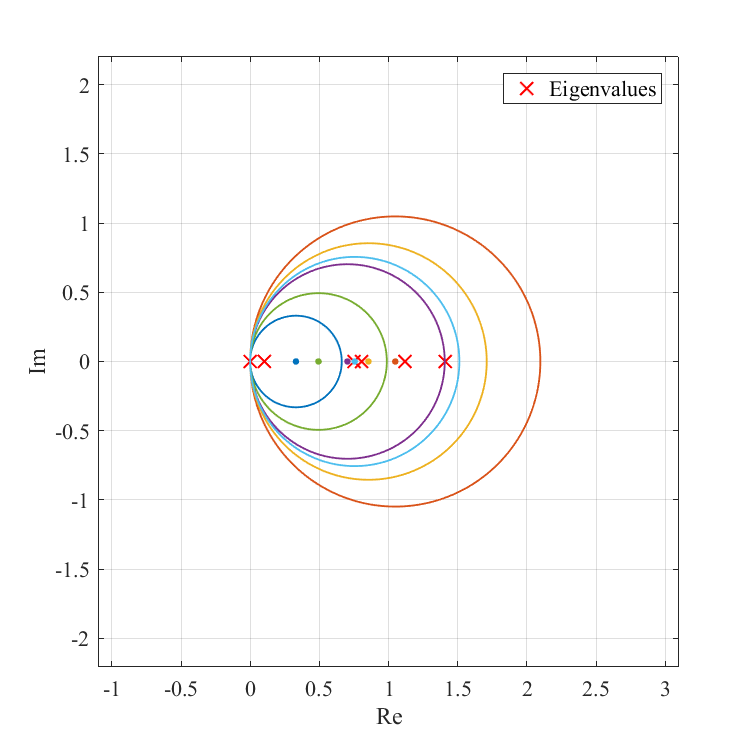}}
	\subfloat[]{
		\includegraphics[trim=0cm 0cm 0cm 0cm, clip=true,  width=4.1cm]{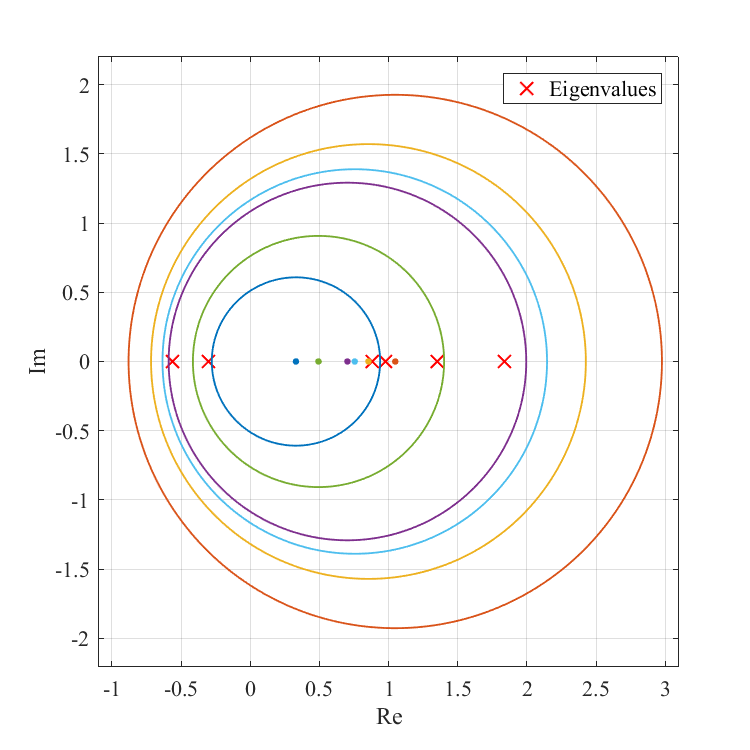}}
	\caption{\small Example~\ref{ex:example1}. Ger\v{s}gorin's disks and eigenvalues of (a): $I-H_1$ (when $\pi=1$). (b): $I-\pi H_1$ (when $\pi>\pi_2$). (c): $L=\Delta-A$ (when $\pi=1$). (d): $\Delta-\pi A$ (when $\pi>\pi_2$).}
	\label{fig:example_n6disks}
\end{figure*}
\end{example}

\begin{example}
\label{ex:example2}
A network of $ n= 20 $ already has $ > 10^6 $ orthants, all potentially containing equilibria of the system \eqref{eqn:SystemSFH}. 
Since, as shown in Example~\ref{ex:example1}, multiple equilibria can appear in the same orthant, this size is already by far out of reach of exhaustive analysis.
The results shown in Fig.~\ref{fig:example_n20} are for a single (nonnegative, irreducible, symmetrizable) realization $ A $, with edges chosen as in an Erd\H{o}s-R\'enyi graph (edge probability $ p=0.1 $) and weights drawn from a uniform distribution. 
All $\psi_i(x_i)$ are chosen equal (again Boltzmann functions).
The results appear to be robust across different realizations of $A$.
Panel (a) of Fig.~\ref{fig:example_n20} shows the number of equilibria (and, in red, the number of orthants to which these equilibria belong) for 500 choices of $ \pi $ uniformly distributed between 1 and 20. For $ \pi $ very small no equilibrium appears, as expected. Equilibria start to appear for values of $ \pi $ that satisfy the condition of Theorem~\ref{theorem:SFnecessaryandSufficientCondition}. 
When $ \pi $ is increased further, then the number of equilibria rapidly grows. For each value of $ \pi $, $ 10^4 $ different initial conditions were tested (we used the {\tt{fsolve}} function of Matlab to compute equilibria). 
The number of equilibria found in this way oscillated between 500 and 600 for a broad range of $ \pi $ values, belonging to 400-500 different orthants. 
As shown in Theorem~\ref{theorem:NormBound}, for each $\pi$ all of the equilibria have a norm which is less than the norm of the corresponding positive/negative equilibrium, see Fig.~\ref{fig:example_n20}(b), where the ratio $ \frac{\| \bar{x} \| } {\| x^+ \| } $ is shown. 
These equilibria were tested for local stability. As shown in panel (c), most but not all of them are unstable, with up to 7 unstable eigenvalues in the Jacobian linearization. 
It is also remarkable that all stable equilibria tend to have high norm $ \| \bar{x} \|  $, i.e., they tend to be near the boundary of the ball of radius $ \| x^+ \|  $ to which they have to belong, see Fig.~\ref{fig:example_n20}(c).
Notice from Fig.~\ref{fig:example_n20}(b) and (c) that equilibria of small norm ratio $ \frac{\| \bar{x} \| } {\| x^+ \| } $, corresponding to small values of $ \pi $, tend also to have an equal ratio of positive and negative components: in Fig.~\ref{fig:example_n20}(c) the radial directions are determined by the fraction of $+ $ and $ - $ signs of an equilibrium, and the bisectrix of the second and fourth quadrant, corresponding to 50\% of $+$ and 50\% of $-$, is where these equilibria tend to be localized.

\begin{figure*}[ht!]
	\centering
	\subfloat[]{
		\includegraphics[trim=0cm 0cm 0cm 0cm, clip=true,  width=5.5cm]{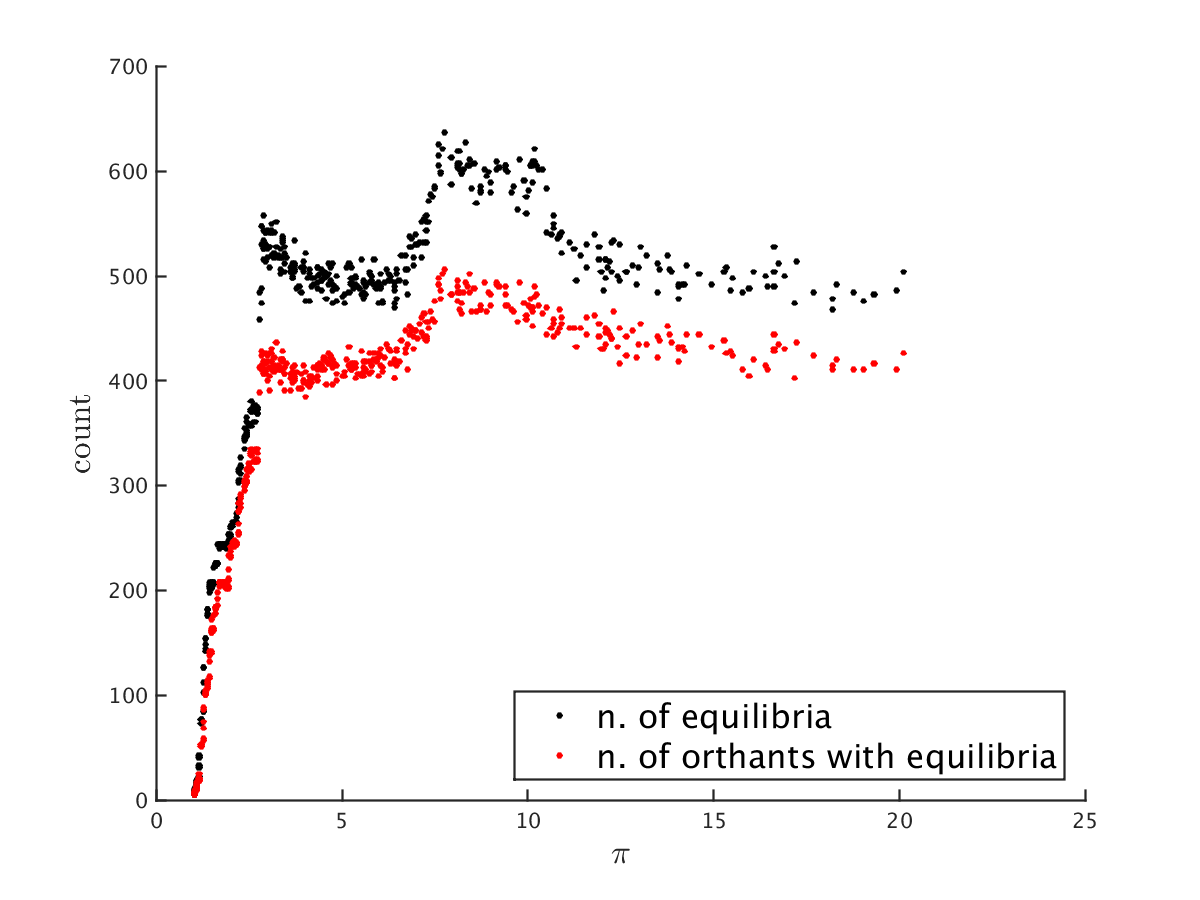}}
	\subfloat[]{
		\includegraphics[trim=0cm 0cm 0cm 0cm, clip=true,  width=5.5cm]{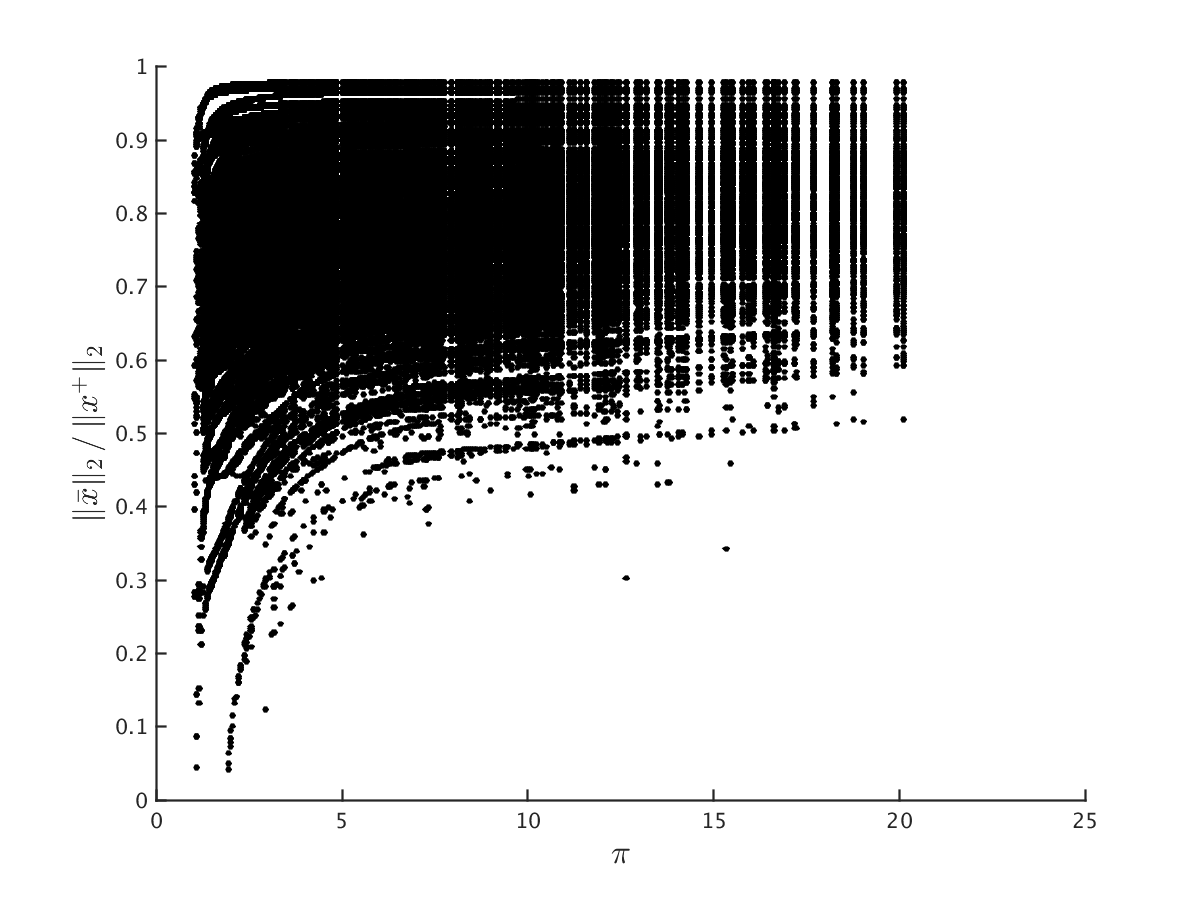}}
	\subfloat[]{
		\includegraphics[trim=0cm 0cm 0cm 0cm, clip=true,  width=5.5cm]{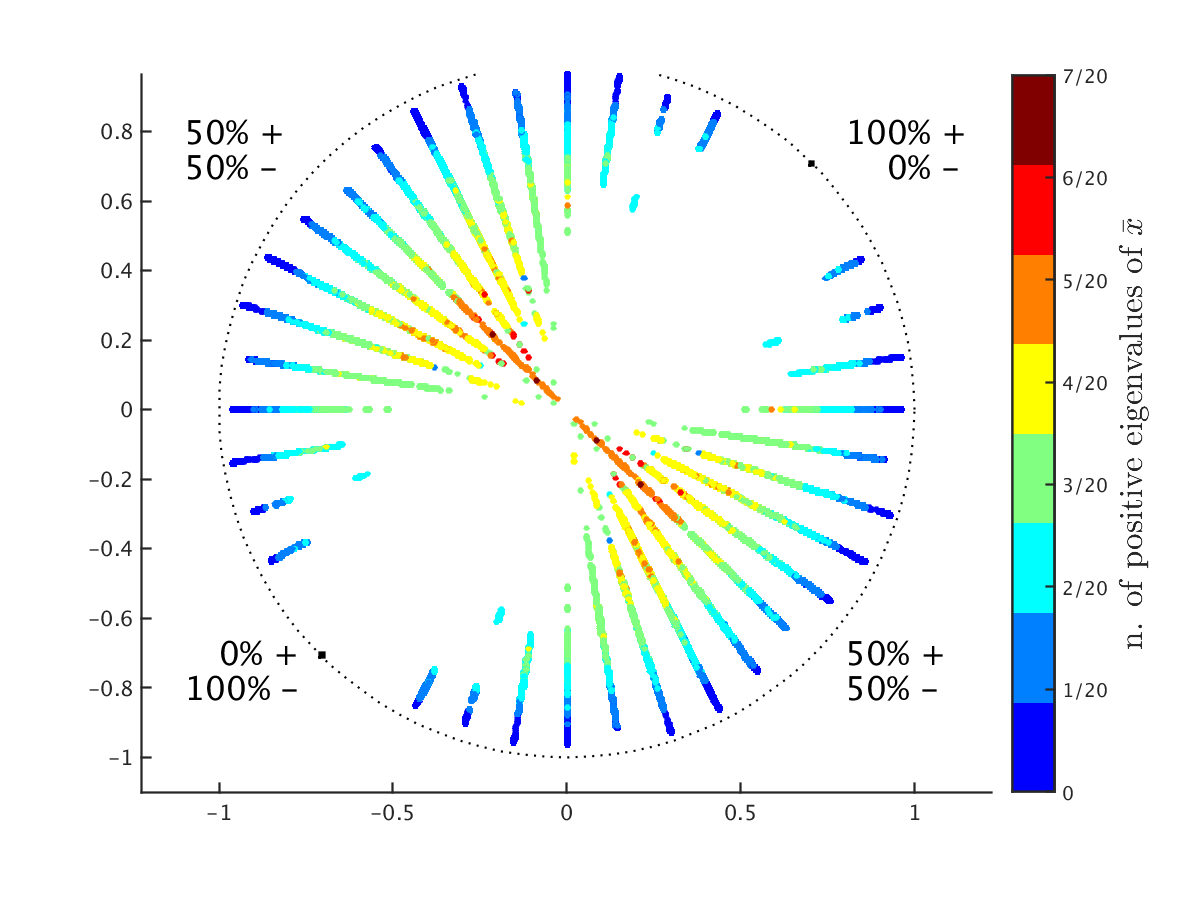}}
	\caption{\small Example~\ref{ex:example2}. Equilibria for system \eqref{eqn:SystemSFA} of size $ n=20$.
	(a): N. of equilibria (black) found on $ 10^4 $ trials for 500 different values of $ \pi$, and n. of orthants to which these equilibria belong (red).
	(b): norm ratio $ \frac{\| \bar{x} \| } {\| x^+ \| } $ of the equilibria $\bar x$ with respect to the corresponding value of $ x^+ $, as $\pi$ changes.
	(c): Distribution of the equilibria according to the norm ratio $ \frac{\| \bar{x} \| } {\| x^+ \| } $ (radial distance from origin), to the number of positive eigenvalues of the Jacobian linearization (colormap), and to the n. of components of negative sign (radial angle). Black squares on the unit disk are $ x^+ $ and $x^- $.}
	\label{fig:example_n20}
\end{figure*}
\end{example}

\section{Final considerations and conclusions}

In this work we have investigated the presence of fixed points for a particular class of nonlinear interconnected cooperative systems, where the nonlinearities are (strictly) monotonically increasing and saturated. 
We have proposed necessary and sufficient conditions on the spectrum of the adjacency matrix of the network which are required for the existence of multiple equilibria not contained in $ \mathbb{R}^n_+/\mathbb{R}^n_-$ when the nonlinearities assume both sigmoidal and nonsigmoidal configuration.
The stability properties of these equilibria have also been investigated.
Although we cannot analytically quantify the number of such equilibria, we can locate them in the solid disk whose radius is given by the norm of the positive equilibrium $ x^+$. 


When interpreted in terms of collective decision-making of agent systems, our results can be recapitulated as follows: 
\begin{itemize}
\item For a low value of the social effort parameter $ \pi $ (i.e., for $ \pi < 1 $), the agents are not committed enough to reach an agreement;
\item For values of $ \pi $ between $1$ and $ \frac{1}{\lambda_{n-1} (H_1)} $, the agents have the right dose of commitment to achieve an agreement among two alternative options $ x^+$, $ x^-$;
\item For values of $ \pi $ bigger than $ \frac{1}{\lambda_{n-1} (H_1) }$, the agents start to become overcommitted, which can lead to other possible decisions, depending on the initial conditions of the system. 
All these extra decisions represent disagreement situations, i..e,  they do not belong to $ \mathbb{R}^n_+ $ or $ \mathbb{R}^n_- $. 
\end{itemize}

Future work includes gaining a better understanding of the bifurcation pattern for $ \pi > \frac{1}{\lambda_{n-1} (H_1) }$, and in presence of ``informed agents'' in the sense of \cite{bib:FranciSrivastavaLeonard2015}. 
It is well-known that the algebraic connectivity is strongly influenced by the topology of the network \cite{Donetti2006Optimal}. We expect that similar arguments apply {\em tamquam} to $\lambda_{n-1} (H_1) $. 
What remains to be checked is whether these topological considerations are applicable to concrete examples of collective decision-making.

\bibliographystyle{IEEEtran}

\begin{thebibliography}{10}
\providecommand{\url}[1]{#1}
\csname url@samestyle\endcsname
\providecommand{\newblock}{\relax}
\providecommand{\bibinfo}[2]{#2}
\providecommand{\BIBentrySTDinterwordspacing}{\spaceskip=0pt\relax}
\providecommand{\BIBentryALTinterwordstretchfactor}{4}
\providecommand{\BIBentryALTinterwordspacing}{\spaceskip=\fontdimen2\font plus
\BIBentryALTinterwordstretchfactor\fontdimen3\font minus
  \fontdimen4\font\relax}
\providecommand{\BIBforeignlanguage}[2]{{%
\expandafter\ifx\csname l@#1\endcsname\relax
\typeout{** WARNING: IEEEtran.bst: No hyphenation pattern has been}%
\typeout{** loaded for the language `#1'. Using the pattern for}%
\typeout{** the default language instead.}%
\else
\language=\csname l@#1\endcsname
\fi
#2}}
\providecommand{\BIBdecl}{\relax}
\BIBdecl

\bibitem{bib:FranciSrivastavaLeonard2015}
\BIBentryALTinterwordspacing
A.~Franci, V.~Srivastava, and N.~E. Leonard, ``{A Realization Theory for
  Bio-inspired Collective Decision-Making},'' \emph{arXiv preprint}, 2015.
  [Online]. Available: \url{https://arxiv.org/abs/1503.08526v1}
\BIBentrySTDinterwordspacing

\bibitem{Gray2017Agent}
R.~Gray, A.~Franci, V.~Srivastava, and N.~E. Leonard, ``An agent-based
  framework for bio-inspired value-sensitive decision-making,'' in
  \emph{Proceedings of the IFAC World Congress}, 2017.

\bibitem{Leonard14Multi}
N.~E. Leonard, ``Multi-agent system dynamics: Bifurcation and behavior of
  animal groups,'' \emph{Annual Reviews in Control}, vol.~38, no.~2, pp.
  171--183, 2014.

\bibitem{bib:AltafiniLini2015}
C.~Altafini and G.~Lini, ``{Predictable dynamics of opinion forming for
  networks with antagonistic interactions},'' \emph{IEEE Transactions on
  Automatic Control}, vol.~60, no.~2, pp. 342--357, 2015.

\bibitem{bib:Altafini2012}
C.~Altafini, ``{Dynamics of opinion forming in structurally balanced social
  networks},'' \emph{PLoS ONE}, vol.~7, no.~6, pp. 5876--5881, 2012.

\bibitem{bib:LiChenAihara2006}
C.~Li, L.~Chen, and K.~Aihara, ``{Stability of genetic networks with SUM
  regulatory logic: Lur'e system and LMI approach},'' \emph{IEEE Transactions
  on Circuits and Systems I: Regular Papers}, vol.~53, no.~11, pp. 2451--2458,
  2006.

\bibitem{bib:Haykin1999}
S.~Haykin, ``{Neural networks: A Comprehensive Foundation},'' pp. 409--412,
  1999.

\bibitem{bib:Hopfield1984}
J.~J. Hopfield, ``{Neurons with graded response have collective computational
  properties like those of two-state neurons.}'' \emph{Proc. Natl. Acad. Sci.
  USA}, vol.~81, no.~10, pp. 3088--3092, 1984.

\bibitem{Kaszkurewicz2000Matrix}
E.~Kaszkurewicz and A.~Bhaya, \emph{Matrix diagonal stability in systems and
  computation}.\hskip 1em plus 0.5em minus 0.4em\relax Boston: Birkh{\"a}user,
  2000.

\bibitem{siljak1978large}
D.~Siljak, \emph{Large-Scale Dynamic Systems: Stability and Structure}.\hskip
  1em plus 0.5em minus 0.4em\relax North-Holland, 1978.

\bibitem{bib:HirschSmith2005}
M.~W. Hirsch and H.~Smith, ``{Monotone dynamical systems},'' in \emph{Handbook
  of differential equations: ordinary differential equations}, A.~Canedea,
  P.~Drabek, and A.~Fonda, Eds.\hskip 1em plus 0.5em minus 0.4em\relax Elsevier
  North Holland, Boston Massachusetts, 2005, vol.~2, ch. Monotone D, pp.
  239--358.

\bibitem{bib:Smith1988}
H.~L. Smith, ``{Systems of Ordinary Differential Equations Which Generate an
  Order Preserving Flow. A Survey of Results},'' \emph{SIAM Review}, vol.~30,
  no.~1, pp. 87--113, 1988.

\bibitem{bib:CohenGrossberg1988}
M.~A. Cohen and S.~Grossberg, \emph{{Artificial Neural Networks: Theoretical
  Concepts}}, V.~Vemuri, Ed.\hskip 1em plus 0.5em minus 0.4em\relax IEEE
  Computer Society Press, 1988.

\bibitem{bib:ChengLinShih2006}
C.-Y. Cheng, K.-H. Lin, and C.-W. Shih, ``{Multistability in Recurrent Neural
  Networks},'' \emph{SIAM Journal on Applied Mathematics}, vol.~66, no.~4, pp.
  1301--1320, 2006.

\bibitem{bib:ChengLinShih2007}
C.~Y. Cheng, K.~H. Lin, and C.~W. Shih, ``{Multistability and convergence in
  delayed neural networks},'' \emph{Physica D: Nonlinear Phenomena}, vol. 225,
  no.~1, pp. 61--74, 2007.

\bibitem{bib:ZengZheng2012}
Z.~Zeng and W.~X. Zheng, ``{Multistability of neural networks with time-varying
  delays and concave-convex characteristics},'' \emph{IEEE Transactions on
  Neural Networks and Learning Systems}, vol.~23, no.~2, pp. 293--305, 2012.

\bibitem{bib:LuWangChen2011}
W.~Lu, L.~Wang, and T.~Chen, ``{On attracting basins of multiple equilibria of
  a class of cellular neural networks},'' \emph{IEEE Transactions on Neural
  Networks}, vol.~22, no.~3, pp. 381--394, 2011.

\bibitem{bib:FortiTesi2015}
M.~Forti and A.~Tesi, ``{New conditions for global stability of neural networks
  with application to linear and quadratic programming problems},'' \emph{IEEE
  Transactions on Circuits and Systems I: Fundamental Theory and Applications},
  vol.~42, no.~7, pp. 354--366, 1995.

\bibitem{bib:ZhangWangLiu2014}
H.~Zhang, Z.~Wang, and D.~Liu, ``{A comprehensive review of stability analysis
  of continuous-time recurrent neural networks},'' \emph{IEEE Transactions on
  Neural Networks and Learning Systems}, vol.~25, no.~7, pp. 1229--1262, 2014.

\bibitem{UgoAbara2016Spectral}
P.~{Ugo Abara}, F.~Ticozzi, and C.~Altafini, ``Spectral conditions for
  stability and stabilization of positive equilibria for a class of nonlinear
  cooperative systems,'' \emph{Automatic Control, IEEE Transactions on}, to
  appear, 2017.

\bibitem{UgoAbara2015Existence}
------, ``Existence, uniqueness and stability properties of positive equilibria
  for a class of nonlinear cooperative systems,'' in \emph{54th IEEE Conference
  on Decision and Control}, Osaka, Japan, 2015.

\bibitem{Deabreu2007Old}
N.~M.~M. de~Abreu, ``Old and new results on algebraic connectivity of graphs,''
  \emph{Linear Algebra and its Applications}, vol. 423, no.~1, pp. 53 -- 73,
  2007.

\bibitem{Donetti2006Optimal}
L.~Donetti, F.~Neri, and M.~A. Muñoz, ``Optimal network topologies: expanders,
  cages, ramanujan graphs, entangled networks and all that,'' \emph{Journal of
  Statistical Mechanics: Theory and Experiment}, vol. 2006, no.~08, p. P08007,
  2006.

\bibitem{Olfati2003Consensus}
R.~Olfati-Saber and R.~Murray, ``Consensus problems in networks of agents with
  switching topology and time-delays,'' \emph{Automatic Control, IEEE
  Transactions on}, vol.~49, no.~9, pp. 1520 -- 1533, sept. 2004.

\bibitem{bib:HornJohnson1990}
R.~A. Horn and C.~R. Johnson, \emph{{Matrix Analysis}}, 1st~ed., C.~U. Press,
  Ed.\hskip 1em plus 0.5em minus 0.4em\relax Cambridge University Press, 1990.

\bibitem{bib:BermanPlemmons1994}
A.~Berman and R.~J. Plemmons, \emph{{Nonnegative matrices in the mathematical
  sciences}}.\hskip 1em plus 0.5em minus 0.4em\relax SIAM, 1994, vol.~9, no.~3.

\bibitem{bib:Gruber1995}
B.~Gruber, \emph{{Symmetries in Science VIII}}, B.~Gruber, Ed.\hskip 1em plus
  0.5em minus 0.4em\relax Springer US, 1995.

\bibitem{bib:DiasCastonguayDourado2016}
\BIBentryALTinterwordspacing
E.~S. Dias, D.~Castonguay, and M.~C. Dourado, ``{Algorithms and Properties for
  Positive Symmetrizable Matrices},'' \emph{TEMA (S{\~{a}}o Carlos)}, vol.~17,
  no.~2, p. 187, 2016. [Online]. Available:
  \url{https://tema.sbmac.org.br/tema/article/view/868}
\BIBentrySTDinterwordspacing

\bibitem{Poignard2017Spectra}
C.~Poignard, T.~Pereira, and J.~P. Pade, ``{Spectra of Laplacian Matrices of
  Weighted Graphs: Structural Genericity Propoerties},''
  \emph{arxiv.org:1704.01677}, 2017.

\bibitem{golubitsky2000singularities}
M.~Golubitsky, I.~Stewart, and D.~Schaeffer, \emph{Singularities and Groups in
  Bifurcation Theory}, ser. Applied Mathematical Sciences.\hskip 1em plus 0.5em
  minus 0.4em\relax Springer New York, 2000, no. v. 1.

\end{thebibliography}

\end{document}